\documentclass[12 pt]{amsart}
\usepackage{amscd,amsmath,latexsym,amssymb}
\usepackage[mathscr]{euscript}
\usepackage[all]{xy}
\usepackage{layout}
\usepackage{pb-diagram}
\usepackage{pstricks}

\newtheorem{theorem}{Theorem}[section]
\newtheorem{corollary}[theorem]{Corollary}
\newtheorem{lemma}[theorem]{Lemma}
\newtheorem{proposition}[theorem]{Proposition}
\theoremstyle{definition}
\newtheorem{definition}[theorem]{Definition}
\theoremstyle{definition}
\newtheorem{example}[theorem]{Example}
\theoremstyle{remark}
\newtheorem*{remark}{Remark}

\def\Proof{\medskip\noindent{\bf Proof: }}
 

\def\Z{\mathbb{Z}}

\def\C{\mathbb{C}}

\def\Q{\mathbb{Q}}
\def\R{\mathbb{R}}

\def\S{\mathbb{S}}

\def\k{\mathfrak{k}}
\def\g{\mathfrak{g}}
\def\t{\mathfrak{t}}
\def\M{\mathcal{M}}

\def\AC{\mathcal{AC}}
\def\ux{\underline{x}}
\def\uy{\underline{y}}

\def\x{\times}

\def\Hom{\text{Hom}}
\def\Ker{\text{Ker}}
\def\Rep{\text{Rep}}
\def\Im{\text{Im}}

\begin{document}

\title[Stable splittings, spaces of representations 
and almost commuting elements]
{Stable splittings, spaces of representations 
and almost commuting elements in Lie groups}

\author[A.~Adem]{Alejandro Adem$^{*}$}
\address{Mathematics Department,
University of British Columbia, Vancouver BC V6T1Z2, Canada}
\email{adem@math.ubc.ca}
\thanks{$^{*}$Partially supported by NSERC}

\author[F.~R.~Cohen]{Frederick R. Cohen$^{**}$}
\address{Department of Mathematics,
University of Rochester, Rochester NY 14627, USA}
\email{cohf@math.rochester.edu}
\thanks{$^{**}$Partially supported by the NSF}

\author[J.~M.~G\'omez]{Jos\'e Manuel G\'omez}
\address{Mathematics Department,
University of British Columbia, Vancouver BC V6T1Z2, Canada}
\email{josmago@math.ubc.ca}

\begin{abstract}
In this paper the space of almost commuting elements in a Lie group is 
studied through a homotopical point of view. In particular a stable 
splitting after one suspension is derived for these spaces and their 
quotients under conjugation. A complete description for the stable
factors appearing in this splitting is provided for compact connected 
Lie groups of rank one.By using symmetric products, the colimits 
$\Rep(\Z^n, SU)$, $\Rep(\Z^n,U)$ and $\Rep(\Z^n, Sp)$ 
are explicitly described as 
finite products of
Eilenberg-MacLane spaces.
\end{abstract}
              
\maketitle

\section{Introduction}
Let $G$ denote a Lie group, for each integer $n\ge 1$, the set 
$\Hom(\Z^{n},G)$ can be topologized in a natural way as a 
subspace of $G^{n}$ and corresponds to the space of ordered 
commuting $n$-tuples in $G$. The group $G$ acts by conjugation on 
$\Hom(\Z^{n},G)$ and the orbit space 
$\Rep(\Z^{n},G):=\Hom(\Z^{n},G)/G$ can be identified with
the moduli 
space of isomorphism classes of flat connections on principal 
$G$-bundles over the torus $(\S^{1})^{n}$.  
Now if $H=G/K$, where $K\subset G$ is a closed subgroup contained
in the center of $G$,
then a commuting  $n$-tuple $(x_{1},...,x_{n})$ of elements in $H$,
lifts to a sequence
$(\tilde{x}_{1},...,\tilde{x}_{n})$ of
elements in $G$ which do not necessarily commute;
instead they form a \textsl{$K$--almost commuting sequence} in 
$G$ in the sense that $[\tilde{x}_{i},\tilde{x}_{j}]\in K$
for all $i$ and $j$. 

Denote the space of all $K$-almost commuting 
$n$-tuples in $G$ by $B_{n}(G,K)$; note that it fits into
a principal $K^n$--bundle
$B_n(G,K) \to \Hom (\Z^n, G/K)$. 
The group $G$ 
acts by conjugation on $B_{n}(G,K)$ with orbit space 
denoted by 
$\bar{B}_{n}(G,K)=B_{n}(G,K)/G$. 
From the standard classification theorems 
(see \cite[Theorem 6.19]{HM}), 
it is known that if $G$ is a compact, connected Lie group, then 
$G=\tilde{G}/K$, where $K$ is 
a finite subgroup in the center of $\tilde{G}$
and $\tilde{G}
=(\S^1)^r\times G_1\times\dots\times G_s$
where the $G_1,\dots , G_s$ are simply connected
simple Lie groups. Thus the $K$--almost commuting
elements clearly play a key role in the analysis of
spaces of
commuting elements for compact Lie groups. This approach
was applied in \cite{BFM} for the more geometrically
accessible situation of commuting triples.
This paper is concerned with the study of the spaces 
of the form $B_{n}(G,K)$ and their orbit spaces 
$\bar{B}_{n}(G,K)$ from
a \textsl{homotopical} point of view for all values of $n\ge 1$.
In particular a basic goal is to compute 
the number of path-connected components and 
determine their stable homotopy types.

\medskip

A basic ingredient is the existence
of a natural simplicial structure on $B_*(G,K)$ 
which provides a filtration for each $B_n(G,K)$ given 
by the degeneracy maps
\[
\{(1_{G},...,1_{G})\}=S_{n}^{n}(G,K)\subset 
S_{n}^{n-1}(G,K)\subset\cdots \subset S_{n}^{0}(G,K)=B_{n}(G,K),
\] 
where $S_{n}^{r}(G,K)$ is the subspace of $B_{n}(G,K)$ of sequences 
$(x_{1},...,x_{n})$ where at least $r$ of the $x_{i}$ equal $1_{G}$.
It turns out that when $G$ is compact the pair 
$(S_{n}^{r}(G,K),S_{n}^{r+1}(G,K))$ is a $G$--equivariant
NDR pair for $0\le r\le n$. 
This technical condition implies (see \cite{ABBCG}, \cite{AC}) 
that the previous filtration splits after one suspension and 
thus the following theorem is obtained

\begin{theorem}\label{introdecomposition}
Suppose that $G$ is a compact Lie group and $K\subset Z(G)$ is 
a closed subgroup. For all $n\ge 0$ there is a natural 
$G$-equivariant homotopy equivalence
\[
\Sigma(B_{n}(G,K))\simeq \bigvee_{1\le r\le n}\Sigma
\left(\bigvee^{\binom{{n}}{{r}}}B_{r}(G,K)/S^{1}_{r}(G,K) \right).
\]
In particular, there is a natural homotopy equivalence
\[
\Sigma(\bar{B}_{n}(G,K))\simeq \bigvee_{1\le r\le n}\Sigma
\left(\bigvee^{\binom{{n}}{{r}}}\bar{B}_{r}(G,K)/
\bar{S}^{1}_{r}(G,K) \right).
\]
\end{theorem}

\medskip

A particularly important case of Theorem \ref{introdecomposition} 
occurs when $K=\{1_{G}\}$. In this case, 
the theorem agrees with \cite[Theorem 1.6]{AC} for compact Lie 
groups with the important additional property
that the splitting map is a $G$-equivariant homotopy equivalence. 
This defines a natural 
homotopy equivalence between the quotient spaces modulo conjugation,
thus yielding a new decomposition
at the level of spaces of representations:

\begin{theorem}
For $G$ a compact Lie group, there is a 
homotopy equivalence
\[
\Sigma\Rep(\Z^{n},G)\simeq \bigvee_{1\le r\le n}
\Sigma\left(\bigvee^{\binom{{n}}{{r}}}\Rep(\Z^{n},G)/
\bar{S}^{1}_{r}(G) \right).
\]
\end{theorem}
When $\Rep(\Z^{n},G)$ is path-connected for all $n\ge 1$ then the
following theorem completely describes this decomposition:

\begin{theorem}
Let $G$ be a compact, connected Lie group such that 
$\Rep(\Z^{n},G)$ is connected for every $n\ge 1$. Let $T$ be a maximal 
torus of $G$ and $W$ the Weyl group associated to $T$. Then 
\[
\Rep(\Z^{n},G)\cong T^{n}/W
\]
with $W$ acting diagonally on $T^{n}$. Moreover there is a
homotopy equivalence
\[
\Sigma \Rep(\Z^{n},G) \simeq \bigvee_{1\le r\le n}
\Sigma(\bigvee^{\binom{{n}}{{r}}} T^{\wedge r}/W),
\]
where $T^{\wedge r}$ is the $r$-fold smash product of $T$ 
and $W$ acts diagonally on $T^{\wedge r}$.
\end{theorem}

The previous theorem applies in particular to the cases  
$G=U(m)$, $Sp(m)$ and $SU(m)$. It turns out that for such 
$G$ the representation spaces $\Rep(\Z^{n},G)$ can be identified 
in terms of symmetric products.

\begin{proposition}\label{introsymmetric}
There are  homeomorphisms
\begin{align*}
\Rep(\Z^{n},U(m))&\cong SP^{m}((\S^{1})^{n}),\\
\Rep(\Z^{n},Sp(m))&\cong SP^{m}((\S^{1})^{n}/\Z/2),
\end{align*}
where $\Z/2$ acts diagonally by complex conjugation on
each factor. Moreover for each $m,n\ge 1$ the determinant
map defines a locally trivial bundle
$$\Rep(\Z^n, SU(m))\to \Rep(\Z^n, U(m))\to (\S^1)^n.$$
\end{proposition}

This identification provides a number of interesting consequences, 
for example 
\begin{corollary}
For every $m\ge 1$ there are homeomorphisms
\[
\Rep(\Z^{2},Sp(m))\cong \C P^{m},
\]
\[
\Rep(\Z^{2},SU(m))\cong \C P^{m-1}.
\]
\end{corollary}
 
Another consequence of Proposition \ref{introsymmetric} is the
determination 
of the homotopy type of the inductive limits of the 
spaces of commuting elements modulo conjugation
for the (special) unitary 
groups and the symplectic groups. Note that this result
for the unitary groups is also 
proved in the recent preprint \cite{Ramras}.

\begin{theorem}
For every $n\ge 1$ there are homotopy equivalences
\begin{align*}
\Rep(\Z^{n},SU)\simeq &\prod_{2\le i\le n}
K(\Z^{\binom{{n}}{{i}}},i),\\
\Rep(\Z^{n},U)\simeq &\prod_{1\le i\le n}
K(\Z^{\binom{{n}}{{i}}},i),\\
\Rep(\Z^{n},Sp)\simeq &\prod_{1\le i \le 
\left\lfloor n/2\right\rfloor}
K(\Z^{\binom{{n}}{{2i}}}\oplus (\Z/2)^{r(2i)},2i),
\end{align*}
where $r(i)$ is an integer given by
\begin{equation*}
r(i)=\left\{
\begin{array}{ccc}
\binom{{n}}{{0}}+\binom{{n}}{{1}}+\cdots+\binom{{n}}{{n-i-1}}& 
\text{ if } &1\le i \le n \\
0& \text{ if } &i> n.
\end{array}%
\right. 
\end{equation*}
\end{theorem}

Also, explicit computations for the stable factors that appear in 
Theorem \ref{introdecomposition} are provided for non-abelian
compact connected rank
one Lie groups. If $G$ is such a group then it is isomorphic to 
$SU(2)$ or $SO(3)$. Let $\lambda_{n}$ denote the canonical vector 
bundle  over the projective space $\R P^{n}$ and for a vector 
bundle $\mu$ over $\R P^{n}$, $(\R P^{n})^{\mu}$ denotes the 
corresponding Thom space. In these cases the factors 
$B_{n}(G,K)/S_{n}^{1}(G,K)$ are given as follows. 

\medskip

\begin{equation*}
\Hom(\Z^{n},SO(3))/S^{1}_{n}(SO(3))\cong \left\{ 
\begin{array}{cl}
\R P^{3} & \text{if }n=1, \\  
(\R P^{2})^{n\lambda_{2}}\vee 
\left(\bigvee_{C(n)} (\S^{3}/Q_{8})_{+}\right)& \text{if }n\ge 2, 
\end{array}%
\right. 
\end{equation*}
where $Q_{8}$ is the quaternion group of order 8 and
\[
C(n)=\frac{1}{2}(3^{n-1}-1).
\]

\medskip

For $SU(2)$, there are two possibilities as $Z(SU(2))\cong \Z/2$: 
\begin{equation*}
\Hom(\Z^{n},SU(2))/S^{1}_{n}(SU(2))\cong\left\{ 
\begin{array}{cl}
\S^{3} & \text{if }n=1, \\  
(\R P^{2})^{n\lambda_{2}}/s_{n}(\R P^{2})& \text{if }n\ge 2. 
\end{array}%
\right. 
\end{equation*}

Here $s_{n}$ is the zero section of the vector bundle $n\lambda_{2}$. 
The other case is the space of almost commuting elements in $SU(2)$:

\begin{equation*}
B_{n}(SU(2),\Z/2)/S^{1}_{n}(SU(2),\Z/2)\cong\left\{
\begin{array}{cl}
\S^{3} & \text{if }n=1, \\  
(\bigvee_{K(n)} \R P^{3}_{+})\vee 
(\R P^{2})^{n\lambda_{2}}/s_{n}(\R P^{2})  & \text{if }n\ge 2, \\  
\end{array}
\right.
\end{equation*}
where
\[
K(n)=\frac{7^{n}}{24}-\frac{3^{n}}{8}+\frac{1}{12}.
\]

\medskip

Finally, using the work in \cite{BFM}, explicit descriptions for 
the stable factors modulo conjugation $\bar{B}_{n}(G,K)/
\bar{S}_{n}^{1}(G,K)$ that appear in Theorem \ref{introdecomposition} 
are given when $G$ is a  compact, connected, simply connected, 
simple Lie group and for $n=1$, $2$ and $3$.

\section{The simplicial spaces $B_{n}(q,G,K)$}

A family of simplicial spaces $B_{n}(q,G,K)$ is defined in 
this section, where $G$ is a Lie group and $K\subset Z(G)$ is  
a closed subgroup. These spaces form a generalization of  
the simplicial spaces $B_{n}(q,G)$ defined in  \cite{ACT} as 
the case $K=\{1_{G}\}$ corresponds to $B_{n}(q,G)$. 
Some general properties are derived.

\medskip

The definition of the spaces $B_{n}(q,G)$ is given first. 
Let $F_{n}$ be the free group on $n$-letters. For each $q\ge 1$, 
let $\Gamma^{q}$ denote the $q$-stage of the descending central 
series for $F_{n}$. The groups $\Gamma^{q}$ are defined inductively 
by setting $\Gamma^{1}=F_{n}$ and $\Gamma^{q+1}=[\Gamma^{q},F_{n}]$. 
As in \cite{ACT} define
\[
B_{n}(q,G):=\Hom(F_{n}/\Gamma^{q},G).
\]
The topology on $B_{n}(q,G)$ is given as follows. An element in 
$\Hom(F_{n}/\Gamma^{q},G)$ can be seen as a homomorphism 
$f:F_{n}\to G$ such that $f(\Gamma^{q})=\{1_{G}\}$ so that 
$f$ descends to a homomorphism $\bar{f}:F_{n}/\Gamma^{q}\to G$. 
The map $f$ is determined by the image of a set of generators 
$e_{1},...,e_{n}$; 
that is, $f$ is determined by the $n$-tuple $(f(e_{1}),...,f(e_{n}))$. 
This gives an inclusion of $\Hom(F_{n}/\Gamma^{q},G)$ into $G^{n}$ 
and via this inclusion $B_{n}(q,G)$ is given the subspace topology. 
The different $B_{n}(q,G)$ form simplicial 
spaces and can be used to obtain a filtration
\[
B(2,G)\subset B(3,G)\subset\cdots \subset B(q,G)
\subset\cdots B(\infty,G)=BG,
\]
of $BG$, where $B(q,G)$ denotes the geometric realization of 
$B_{n}(q,G)$. 

\begin{definition}
Let $Q$ be a group. A sequence of subgroups $\Gamma^{r}(Q)$ 
is defined inductively as follows. For $r=1$, $\Gamma^{1}(Q)=Q$ 
and $\Gamma^{i+1}(Q)=[\Gamma^{i}(Q),Q]$. The different 
$\Gamma^{r}(Q)$ form the descending central series of $Q$
\[
\cdots \subset \Gamma^{i+1}(Q)\subset\Gamma^{i}(Q)
\subset \cdots \Gamma^{2}(Q)\subset \Gamma^{1}(Q)\subset Q.
\]
\end{definition}

A discrete group $Q$ is said to be nilpotent if 
$\Gamma^{m+1}(Q)=\{1\}$ for some $m$ and the least such $m$ is called 
the nilpotency class of $Q$. Suppose that $\bar{f}:F_{n}/\Gamma^{q}
\to G$ is a homomorphism induced by $f:F_{n}\to G$. If $Q=\Im f$, 
then $f(\Gamma^{q})=\Gamma^{q}(Q)=\{1_{G}\}$; that is, $Q$ is 
nilpotent of nilpotency class less than $q$.  Note that if $q\ge 2$ 
and if $G$ is a topological group of nilpotency less class 
than $q$, then $B(q,G)=BG$. 

\medskip

Suppose now that $H$ is a Lie group of the form $G/K$, for a 
closed subgroup $K\subset Z(G)$; that is, suppose that 
$H$ fits into a central extension of Lie groups
\[
1\to K\to G\to H\to 1.
\]
Take $\bar{f}:F_{n}/\Gamma^{q}\to H$ an element in $B_{n}(q,H)$ 
induced by a homomorphism $f:F_{n}\to H$. Then $f$ is determined 
by the $n$-tuple $(x_{1},...,x_{n})=(f(e_{1}),...,f(e_{n}))$ of 
elements in $H$. Since $K\subset Z(G)$ is a closed subgroup, the 
natural map $p:G\to G/K$ is both a homomorphism and a principal 
$K$-bundle. For each $1\le i\le n$, take $\tilde{x_{i}}$ a lifting 
of $x_{i}$ in $G$. The elements $\tilde{x}_{1},...,\tilde{x}_{n}$ 
define a map $g:F_{n}\to G$. Note however that this map does not 
necessarily descend to a map $F_{n}/\Gamma^{q}\to G$. Instead,  
$p(g(\Gamma^{q}))=f(\Gamma^{q})=\{1\}\in G/K$ and it follows that 
$g(\Gamma^{q})\subset K$. This motivates the following definition.

\begin{definition}
Let $G$ be a Lie group, $K\subset Z(G)$ a closed subgroup and 
$n,q\ge 1$. Define
\[
B_{n}(q,G,K)=\{f:F_{n}\to G ~|~ f \text{ is a homomorphism and } 
f(\Gamma^{q})\subset K\}.
\] 
\end{definition}

As in the case of $B_{n}(q,G)$, the set $B_{n}(q,G,K)$ is topologized 
as a subspace of $G^{n}$. Note that if $K=\{1_{G}\}$ is the trivial 
subgroup, then $B_{n}(q,G)=B_{n}(q,G,\{1_{G}\})$ and  when $q=2$, 
the space $B_{n}(2,G,K)$ consists of the $n$-tuples 
$(x_{1},...,x_{n})$ in $G$ such that $[x_{i},x_{j}]\in K$ for 
all $i$ and $j$. 

\begin{lemma}\label{covering}
Suppose that $G$ is a Lie group, $K\subset Z(G)$ a 
closed subgroup and $q\ge 1$. Then the restriction of the natural 
map $p^{n}:G^{n}\to (G/K)^{n}$ to $B_{n}(q,G,K)$ defines a 
$G$-equivariant $\Hom(F_{n}/\Gamma^{q},K)\cong K^{n}$-principal 
bundle
\[
\phi_{n}:B_{n}(q,G,K)\to B_{n}(q,G/K).
\]  
\end{lemma}
\Proof
The natural map $p:G\to G/K$ is a principal $K$-bundle 
and since $K$ is central, the conjugation action of $G$ 
induces an action of $G$ on $G/K$ such that $p$ is $G$-equivariant. 
It follows that $p^{n}:G^{n}\to (G/K)^{n}$ is a 
$G$-equivariant principal $K^{n}$-bundle. By restriction, the map 
\[
\phi_{n}=(p^{n})_{|(p^{n})^{-1}(B_{n}(q,G/K))}
:(p^{n})^{-1}(B_{n}(q,G,K))\to B_{n}(q,G/K)
\] 
is also a principal $K^{n}$-bundle. 
Thus it is only necessary to show that 
\[ 
(p^{n})^{-1}(B_{n}(q,G/K))=B_{n}(q,G,K).
\]
From the definition it is clear that 
$B_{n}(q,G,K)\subset (p^{n})^{-1}(B_{n}(q,G/K))$. 
On the other hand, take $(x_{1},...,x_{n})\in G^{n}$ such that 
$(p(x_{1}),...,p(x_{n}))\in B_{n}(q,G)$. This is equivalent 
to saying that the subgroup $\bar{Q}\subset G/K$ generated by 
$p(x_{1}),...,p(x_{n})$ has nilpotency class less than $q$. 
If $Q$ if the subgroup of $G$ generated by $x_{1},...,x_{n}$, 
then $p(\Gamma^{q}Q)=\Gamma^{q}(p(Q))=\Gamma^{q}(\bar{Q})=
\{1_{G}\}$. It follows that $\Gamma^{q}Q\subset K$ and thus 
$(x_{1},...,x_{n})\in B_{n}(q,G,K).$  
\qed
\medskip

\begin{remark} In \cite[Lemma 2.2]{Goldman}, Goldman showed that 
if $\pi$ is a finitely generated group and $p:G'\to G$ is a 
local isomorphism, then  composition with $p$ defines a continuous 
map 
\[
p_{*}:\Hom(\pi,G')\to \Hom(\pi,G),
\]
such that the image of $p_{*}$ is a union of connected components 
of $\Hom(\pi,G)$ and if $Q$ is a connected component in the image 
of $p_{*}$, then the restriction 
\[
(p_{*})_{|p_{*}^{-1}(Q)}:p_{*}^{-1}(Q)\to Q
\] 
is a covering space with covering group $\Hom(\pi,K)$, where 
$K=\text{Ker}(p)$. In particular, this 
applies for $\pi =F_{n}/\Gamma_{q}$. This shows that in this case 
$p$ defines a continuous map 
\[
p_{*}:B_{n}(q,G')\to B_{n}(q,G),
\]
such that the image of $p_{*}$ is a union of connected components 
of $B_{n}(q,G)$. However, the map $p_{*}$ is not surjective in 
general. For example, if $G'=SU(2)$ and $G=G'/K$, where 
$K=Z(G)=\{\pm 1\}$, then $G'/K=SU(2)/\{\pm 1\}\cong SO(3)$, 
$B_{2}(2,SU(2))=\Hom(\Z^{2},SU(2))$ is path-connected but 
$B_{2}(2,SO(3))=\Hom(\Z^{2},SO(3))$ has two path-connected 
components. The lack of surjectivity of $p_{*}$ is precisely what 
motivates the study of the spaces $B_{n}(q,G,K)$. 
\end{remark}

Take now $G$ a Lie group and $K\subset Z(G)$ a closed subgroup. 
For each $q\ge 2$ the collection $B_{*}(q,G,K)$ can be seen as a 
simplicial space. Indeed, recall that for any Lie group $G$ there 
is an associated simplicial space $B_{*}G$, where $B_{n}G=G^{n}$, 
\[
s_{j}(g_{1},...,g_{n})=(g_{1},...,g_{j},1_{G},g_{j+1},...,g_{n}),
\] 
and
\begin{equation*}
\partial _{i}(g_{1},...g_{n})=\left\{
\begin{array}{cl}
(g_{2},...,g_{n}) & \text{ if }i=0, \\
(g_{1},...,g_{i}g_{i+1},...,g_{n}) & \text{ if }0<i<n, \\
(g_{1},...,,g_{n-1}) & \text{ if }i=n.
\end{array}
\right.
\end{equation*}
The restriction of the maps $s_{j}$ and $\partial_{i}$ 
keep the different $B_{n}(q,G,K)$ invariant as the maps 
$s_{j}$ and $\partial_{i}$ are induced by homomorphisms between 
free groups and if $f:A\to B$ is a homomorphism of groups, 
then $f(\Gamma^{q}(A))\subset \Gamma^{q}(B)$. 
Note that the conjugation action of $G$ on each 
$B_{n}(q,G,K)$ makes $B_{*}(q,G,K)$ into a $G$--simplicial 
space. Moreover, the different maps 
\[
\phi_{n}:B_{n}(q,G,K)\to B_{n}(q,G/K)
\]  
define a $G$-equivariant simplicial map.

\medskip 

The spaces of the form $B_{n}(2,G,K)$ are the subject of study 
of this article and from now on only the case $q=2$ will be 
considered. The abbreviated notation $B_n(G,K)$ from the 
introduction will be used from now on.

\section{Almost commuting elements in Lie groups}\label{section 1}

In this section, following the approach in \cite{BFM}, 
a decomposition for the spaces $B_{n}(G,K)$ is given 
by keeping track of the different commutators. This decomposition 
gives a good way of detecting connected components for the space 
$\Hom(\Z^{n},G/K)$.

\medskip

Suppose that $G$ is a Lie group and that 
$K\subset Z(G)$ is a closed subgroup. 
By definition $B_{n}(G,K)$ is the space of homomorphisms 
$f:F_{n}\to G$ such that $f(\Gamma^{2})\subset K$. 
Since $\Gamma^{2}=[F_{n},F_{n}]$ it follows that under 
the identification of $B_{n}(G,K)$ as a subspace of $G^{n}$, 
an $n$-tuple $\ux:=(x_{1},...,x_{n})\in G^{n}$ belongs to 
$B_{n}(G,K)$ if and only if $d_{ij}:=[x_{i},x_{j}]\in K$; that is, 
the elements $x_{i}$ and $x_{j}$ 
commute up to elements in $K$. Such a sequence is called an 
almost commuting sequence.

\begin{definition}
Suppose $G$ is a Lie group and $K\subset Z(K)$ is a subgroup. 
A sequence $(x_{1},...,x_{n})\in G$ is said to be $K$-almost 
commuting if $[x_{i},x_{j}]\in K$ for all $i$ and $j$.
\end{definition}

According to the previous discussion, $B_{n}(G,K)$ is precisely 
the space of $K$-almost commuting $n$-tuples in $G$. Take 
$\ux=(x_{1},...,x_{n})\in B_{n}(G,K)$, therefore 
$d_{ij}=[x_{i},x_{j}]\in K$. The elements $d_{ij}$ are such 
that $d_{ii}=1$ and $d_{ij}=d_{ji}^{-1}$, this means that the 
matrix $D=(d_{ij})$, which is called the type of $\ux$, 
is an antisymmetric matrix with entries in $K\subset Z(G)$. 
Clearly the entries $d_{ij}=d_{ij}(\ux)$ are continuous 
functions of $\ux$ with values in $K$. A decomposition of the 
space $B_{n}(G,K)$ is obtained by keeping track of the different 
commutators $d_{ij}$. More precisely, let $\pi_{0}:K\to \pi_{0}(K)$ 
be the map that identifies the connected components of elements 
in $K$ and define $T(n,\pi_{0}(K))$ to be the set of all $n\x n$ 
antisymmetric matrices $C=(c_{ij})$ with 
entries in $\pi_{0}(K)$. Given $C\in T(n,\pi_{0}(K))$ define 
\[
\AC_{G}(C)=\{(x_{1},...,x_{n})\in G^{n} ~|~ 
\pi_{0}([x_{i},x_{j}])=c_{ij} \};
\]
that is, $\AC_{G}(C)$ is the set of almost commuting sequences in 
$G$ of type $D$ with $\pi_{0}(D)=C$. Each $\AC_{G}(C)$ is a 
subspace of $B_{n}(G,K)$ that is both open and closed because 
$[*,*]$ is a continuous function and $\pi_{0}(K)$ is discrete. 
It follows that each $\AC_{G}(C)$ is a union of connected 
components of $B_{n}(G,K)$. In addition, the conjugation action 
of $G$ on $G^{n}$, restricts to an action of $G$ on each 
$\AC_{G}(C)$. Indeed, if $[x_{i},x_{j}]=d_{ij}\in Z(G)$, 
then $[gx_{i}g^{-1},gx_{j}g^{-1}]=gd_{ij}g^{-1}=d_{ij}$. 
Following \cite{BFM} the orbit space of $\AC_{G}(C)$ under this 
action is denoted by
\[
\M_{G}(C):=\AC_{G}(C)/G.
\]
Note that $(Z(G))^{n}$ acts on each $B_{n}(G,K)$ by left 
componentwise multiplication. Also, when $K=\{1_{G}\}$ the space 
$B_{n}(G,\{1_{G}\})$ is precisely $\Hom(\Z^{n},G)$. 
The orbit space under the action of $G$ is denoted by
\[
\bar{B}_{n}(G,K):=B_{n}(G,K)/G
\]
and when $K=\{1_{G}\}$ 
\[
\Rep(\Z^{n},G):=\Hom(\Z^{n},G)/G.
\]
By definition, the different $\AC_{G}(C)$ give a decomposition 
of $B_{n}(G,K)$ into subspaces that are both open and closed
as follows:
\begin{equation}\label{decomposition for almost}
B_{n}(G,K)=\bigsqcup_{C\in T(n,\pi_{o}(K))} \AC_{G}(C).
\end{equation}
The natural map $p:G\to G/K$ is an open map and the restriction of 
$p^{n}$ to $B_{n}(G,K)$ is a principal $K^{n}$-bundle 
$\phi_{n}:B_{n}(G,K)\to \Hom(\Z^{n},G/K)$ by Lemma \ref{covering}
and is also an open map. This shows that if $C\in T(n,\pi_{0}(K))$, 
then  
\[
\Hom(\Z^{n},G/K)_{C}:=p^{n}(\AC_{G}(C))
\]
is an open subset of $\Hom(\Z^{n},G/K)$ and 
(\ref{decomposition for almost}) induces a decomposition
\begin{equation}\label{decomposition of Hom}
\Hom(\Z^{n},G/K)=\bigsqcup_{C\in T(n,\pi_{0}(K))}\Hom(\Z^{n},G/K)_{C}.
\end{equation}
The complement of each $\Hom(\Z^{n},G/K)_{C}$ is open since it is a 
union of spaces of the form $\Hom(\Z^{n},G/K)_{C'}$, therefore each 
$\Hom(\Z^{n},G/K)_{C}$ is both open and closed in $\Hom(\Z^{n},G/K)$ 
hence a union of connected components. Moreover, for each element 
$C\in T(n,\pi_{o}(K))$ the natural map 
\[
\AC_{G}(C)\to \Hom(\Z^{n},G/K)_{C}
\]
is a principal $K^{n}$-bundle by Lemma \ref{covering}. The 
decomposition (\ref{decomposition of Hom}) gives a good way of 
detecting connected path components of the spaces of the type 
$\Hom(\Z^{n},G/K)$. 
Let $\mathcal{N}(G,K)$ be the number of different $C\in 
T(n,\pi_{0}(K))$ for which $\AC_{G}(C)$ is nonempty. Then the 
following is an immediate corollary.  

\begin{corollary}
Suppose that $G$ is a Lie group and $K\subset Z(G)$ is a closed 
subgroup. Then $\Hom(\Z^{n},G/K)$ has at least $\mathcal{N}(G,K)$ 
connected components.
\end{corollary}

\begin{example} $SO(m)$ fits into a central extension of Lie groups
\[
1\to \Z/2\stackrel{i}{\rightarrow} Spin(m)\stackrel{p}{\rightarrow} 
SO(m)\to 1.
\]
Composition with $p$ gives an exact sequence of pointed sets
\[
1\to \Hom(\Z^{n},\Z/2)\to \Hom(\Z^{n},Spin(m))
\stackrel{p_{*}}{\rightarrow} \Hom(\Z^{n},SO(m)).
\]
In this case $p_{*}$ is not surjective and the image of $p_*$ does 
not contain all the path-connected components of 
$\Hom(\Z^{n},SO(m))$.
As in (\ref{decomposition of Hom}) 
there is a decomposition 
\begin{equation}\label{the case spin(m)}
\Hom(\Z^{n},SO(m))=\bigsqcup_{C\in T(n,\Z/2)}
\Hom(\Z^{n},SO(m))_{C},
\end{equation}
where each $\Hom(\Z^{n},SO(m))_{C}$ is a union 
of connected components of $\Hom(\Z^{n},SO(m))$.
Given a homomorphism $f:\Z^n\to SO(m)$, it lifts compatibly 
to a central extension 
$$1\to \Z/2\to \Gamma \to \Z^n\to 1$$
which maps to $Spin(m)$. Such an extension is determined
by an element in $H^2(\Z^n,\Z/2)\cong (\Z/2)^{\binom{n}{2}}$.
Fixing a basis $x_1,\dots , x_n$ for $\Z^n$, the
components for this cohomology class
are determined by taking liftings $\tilde{x_i}$
and $\tilde{x_j}$ of $x_i$ and $x_j$ respectively
and calculating their commutator. This is how the
transgression map $H^1(\Z/2,\Z/2)\to H^2(\Z^n,\Z/2)$ associated
to the extension $\Gamma$ is
computed in group cohomology. The values obtained are
thus the same as those used to define the matrix
$C\in T(n,\Z/2)$ associated to the almost commuting $n$--tuple
$\tilde{x_1},\dots , \tilde{x_n}$; moreover,
note that the two sets
$H^2(\Z^n,\Z/2)$ and $T(n,\Z/2)$ have the same cardinality. 
By construction the
$2$--dimensional cohomology class associated to $\Gamma$
must be the pullback of the cohomology class which determines
$Spin(m)$; this is precisely the second Stiefel--Whitney class.
More precisely, given $f:\Z^n\to SO(m)$, the class 
$w_2(f)\in H^2(\Z^n,\Z/2)$ is defined as $f^*(w_2)$;
where $w_2\in H^2(BSO(m),\Z/2)$ is the universal Stiefel--Whitney
class.

What this shows is that the decomposition of $\Hom(\Z^n,SO(m))$ 
described above is the decomposition determined
by the second Stiefel--Whitney class; in other words
$f$ and $g$ are in the same $C$--component if and only
if $w_2(f)=w_2(g)$. Moreover in \cite{ACC} it was proved
that for $m>>0$ all the values of $w_2$ can be attained
by suitable homomorphisms, thus establishing that 
$\Hom(\Z^n,SO(m))$ has at least $\mathcal{N}(Spin(m),\Z/2)=
2^{\binom{n}{2}}$ non--empty components. 
Note that the component corresponding to the trivial cohomology
class is precisely $p_*(\Hom(\Z^n,Spin(m)))$, i.e. those
homomorphisms that can be lifted to $Spin(m)$ (which
happens if and only if $w_2=0$).
\end{example}

\begin{example} In \cite{ACG2} the authors use the decomposition 
(\ref{decomposition of Hom})
to study the spaces of the form $\Hom(\Z^{n},G_{m,p})$, where 
for a prime number $p$ and an integer $m\ge 1$, $G_{m,p}$ denotes 
a central product of $m$-copies of $SU(p)$. There it is shown 
that $\Hom(\Z^{n},G_{m,p})$ has
\[
N(n,m,p)=\frac{p^{(m-1)(n-2)}(p^{n}-1)(p^{n-1}-1)}{p^{2}-1}+1
\]
path-connected components of which $N(n,m,p)-1$ are homeomorphic to
\[
(SU(p))^{m}/((\Z/p)^{m-1}\times E_{p}),
\]
where $E_{p}\subset SU(p)$
is the quaternion group $Q_{8}$ of order eight when $p=2$ and the
extraspecial $p$--group of order $p^3$ and exponent $p$ when $p>2$.
\end{example}

\section{NDR pairs and compact Lie groups}

In this section it is shown that if $G$ is a compact Lie group 
and $K\subset Z(G)$ is a closed subgroup, then the simplicial 
spaces $B_{*}(G,K)$ are simplicially NDR as defined in \cite{ABBCG}.

\medskip

To start the usual definitions of $G$-NDR pairs and 
strong NDR pairs are given.

\begin{definition}
Suppose that $G$ is a topological group and take $(X,A)$ a 
$G$-pair. The pair $(X,A)$ is said to be a $G$-NDR pair if 
there exists a $G$-equivariant continuous function 
\[
h:X\x [0,1]\to X,
\]
with $G$ acting trivially on $[0,1]$ and a $G$-invariant 
continuous map $u:X\to [0,1]$ such that the following 
conditions are satisfied.
\begin{enumerate}
\item $A=u^{-1}(0)$,
\item $h(x,0)=x$ for all $x \in X$,
\item $h(a,t)=a$ for all $a\in A$ and all $t\in [0,1]$, and
\item $h(x,1)\in A$ for all $x\in u^{-1}(\left[0,1\right)).$
\end{enumerate}
\end{definition}

Note that when $G=\{1\}$ is the trivial group then 
the previous definition corresponds to Steenrod's original 
definition of an NDR pair as in \cite{Steenrod}. The pair 
$(h,u)$ is called a $G$-NDR representation for the 
topological pair $(X,A)$. The following is a straightforward 
lemma.

\begin{lemma}\label{G-NDR}
Suppose that $(X,A)$ is a $G$-NDR pair represented by $(h,u)$. 
By passing to orbit spaces, the pair $(h,u)$ induces an NDR 
representation $(\bar{h},\bar{u})$ for the pair $(X/G,A/G)$. 
Also, for any subgroup $H\subset G$, by passing to fixed 
points $(h,u)$ induces an NDR representation $(h^{H},u^{H})$ 
for the pair $(X^{H},A^{H})$.
\end{lemma}

In \cite{MayGeom}, May considered a more restricted version 
of NDR pairs to study simplicial spaces, in particular to 
study the natural filtration of the geometric realization of 
a simplicial space. The precise definition is the following.

\begin{definition}
An NDR pair $(X,A)$ represented by the pair $(h,u)$ is said to 
be a strong NDR pair if $u(h(x,t))<1$ for all $t\in [0,1]$, 
whenever $u(x)<1$. A $G$-NDR pair $(X,A)$ is a strong $G$-NDR pair 
if $(h,u)$ is a $G$-NDR representation of $(X,A)$ 
that is also a strong NDR pair representation for $(X,A)$.
\end{definition}

Given an NDR pair $(X,A)$ represented by $(h,u)$, an integer 
$n\ge 0$ and $1\le r\le n$, denote by $F_{r}X^{n}$ the subspace 
of $X^{n}$ consisting of $n$-tuples $(x_{1},...,x_{n})$ 
where at least $r$ of the $x_{j}$'s are in $A$. In \cite{AC}, 
it is proved that $(X^{n},F_{r}X^{n})$ is an NDR-pair. 
This naturally extends to the equivariant situation as it is 
shown next.

\begin{lemma}\label{AC NDRpairs}
Suppose that $(X,A)$ is a $G$-NDR pair. 
Then the pair $(X^{n},F_{r}X^{n})$ is a 
$G\times\Sigma_{n}$-equivariant NDR pair represented by 
$(h_{n}, u_{r,n})$, where 
\[
u_{r,n}(x_{1},...,x_{n})=\frac{1}{r}\min_{1\le i_{1}<
\cdots< i_{r}\le n} \left(u(x_{i_{1}})+\cdots+u(x_{i_{r}})\right)
\]
and
\[
h_{n}((x_{1},...,x_{n}),t)=(h(x_{1},t_{1}),...,h(x_{n},t_{n})),
\]
with
\begin{equation*}
t_{i}=\left\{ 
\begin{array}{ll}
t\min_{i\neq n}(u(x_{m})/u(x_{i})) & \text{if some }
u(x_{m})<u(x_{i}),\
m\neq i, \\ 
t & \text{if all }u(x_{m})\geq u(x_{i}).%
\end{array}%
\right. 
\end{equation*}
\end{lemma}
\Proof By \cite[Lemma 7.2]{AC} $(h_{n},u_{r,n})$ 
is a $\Sigma_{n}$-NDR-pair representation for 
$(X^{n},F_{r}X^{n})$. Thus it is only necessary to 
see that $(h_{n},u_{r,n})$ respects the $G$-action, where $G$ 
acts diagonally on $X^{n}$. To see this note that $F_{r}X^{n}$ is 
$G$-invariant as $A$ is $G$-invariant, also for any $g\in G$
\begin{align*}
u_{r,n}(g\cdot x_{1},...,g\cdot x_{n})&=\frac{1}{r}
\min_{1\le i_{1}<\cdots< i_{r}\le n}\left(u(g\cdot x_{i_{1}})
+\cdots+u(g\cdot x_{i_{r}})\right)\\
&=\frac{1}{r}\min_{1\le i_{1}<\cdots< i_{r}\le n}
\left(u(x_{i_{1}})+\cdots+u(x_{i_{r}})\right)
\end{align*}  
as $u$ is $G$-invariant. On the other hand, 
for any $t\in [0,1]$ and $t_{i}$ as above for $0\le i\le n$, 
\begin{align*}
h_{n}((g\cdot x_{1},...,g\cdot x_{n}),t)&=
(h(g\cdot x_{1},t_{1}),...,h(g\cdot x_{n},t_{n}))\\
&=(g\cdot h(x_{1},t_{1}),...,g\cdot h(x_{n},t_{n}))\\
&=g\cdot h_{n}((x_{1},...,x_{n}),t).
\end{align*}
Therefore  $(h_{n}, u_{r,n})$ is a $G\times \Sigma_{n}$-NDR 
pair representation for $(X^{n},F_{r}X^{n})$.
\qed
\medskip

The case $X=G$ and $A=\{1_{G}\}$, where $G$ is a Lie group 
is of particular interest. In this case $F_{r}G^{n}\subset G^{n}$ 
is the subspace of $n$-tuples $(x_{1},...,x_{n})$ with 
at least $r$ of the $x_{i}$'s equal to $1_{G}$. Let 
$S_{n}^{r}(G,K):=B_{n}(G,K)\cap F_{r}G^{n}$, thus 
$S_{n}^{r}(G,K)\subset G^{n}$ is the subspace of $K$-almost 
commuting $n$-tuples with at least $r$ components equal to $1_{G}$. 
Note that $G$ acts by conjugation on itself and that each 
$S_{n}^{r}(G,K)$ is invariant under this action. Denote by 
$\bar{S}_{n}^{r}(G,K)$ the orbit space $S_{n}^{r}(G,K)/G$. 
For any Lie group $G$, the pair $(G,\{1_{G}\})$ is an NDR pair 
and thus by the previous lemma $(G^{n},S_{n}^{r}(G))$ is an NDR 
pair. The following proposition provides a criterion that can be 
used to show that each pair $(S_{n}^{r-1}(G,K),S_{n}^{r}(G,K))$ is 
a strong $G$-NDR pair when $G$ is a compact Lie group.

\begin{proposition}\label{NDR centralizers}
Let $G$ be a Lie group and $K\subset Z(G)$ a closed subgroup. 
Suppose $(h,u)$ is a representation of $(G,\{1_{G}\})$ as a strong 
NDR pair with the following additional properties: 
for each $g\in G-\{1_{G}\}$ and any $0\le t<1$ 
\begin{enumerate}
\item $Z(g)=Z(h(g,t))$ and
\item $\{x\in G ~|~ [x,g]\in K\}=\{x\in G ~|~ [x,h(g,t)]\in K\}$.
\end{enumerate}
Then for each $0\le r\le n$ the representation 
$(h_{n},u_{r,n})$ as above restricts to a representation of 
$(S_{n}^{r-1}(G,K),S_{n}^{r}(G,K))$ as a $\Sigma_{n}$-equivariant 
strong NDR pair. Moreover, if the representation $(h,u)$ is 
$G$-equivariant, $G$ acting on itself by conjugation, then the pair 
$(S_{n}^{r-1}(G,K),S_{n}^{r}(G,K))$ is a strong $G$-NDR pair.
\end{proposition} 
\Proof
Suppose first that $(x_{1},...,x_{n})\in  B_{n}(G,K)$; 
that is, $[x_{i},x_{j}]\in K$ for all $i$ and $j$. 
Thus if $t_{i}$ and $t_{j}$ are as explained 
before for each $i$ and $j$ then
\[
x_{i}\in \{x\in G ~|~ [x,x_{j}]\in K\}=\{x\in G ~|~ 
[x,h(x_{j},t_{j})]\in K\},
\]
where the last equality holds by assumption $(2)$. Thus 
\[
[h(x_{j},t_{j}),x_{i}]=[x_{i},h(x_{j},t_{j})]^{-1}\in K
\] 
and therefore
\[
h(x_{j},t_{j})\in \{x\in G ~|~ [x,x_{i}]\in K\}=\{x\in G ~|~ 
[x,h(x_{i},t_{i})]\in K\}.
\] 
From here it can be concluded that 
$[h(x_{i},t_{i}),h(x_{j},t_{j})]\in K$ for all $i$ and $j$, hence 
\[
h_{n}((x_{1},...,x_{n}),t)=(h(x_{1},t_{1}),...,h(x_{n},t_{n}))
\in B_{n}(G,K).
\]
Therefore $h_{n}:B_{n}(G,K)\x [0,1]\to B_{n}(G,K)$ is well 
defined. Moreover each $S_{n}^{r}(G,K)$ is invariant under 
each $h_{n}(\cdot,t)$ as if $(x_{1},...,x_{n})\in S_{n}^{r}(G,K)$, 
then at least $r$ of the $x_{i}$'s are equal to $1_{G}$. Since 
$h(1_{G},t)=1_{G}$ for all $t$, it follows that 
$h_{n}((x_{1},...,x_{n}),t)$ has at least $r$ components equal to 
$1_{G}$. Let $(\tilde{h}_{n},\tilde{u}_{r,n})$ be the restriction 
of $(h_{n},u_{r,n})$ to $S_{n}^{r-1}(G,K)$. Then 
$(\tilde{h}_{n},\tilde{u}_{r,n})$ represents the pair 
$(S_{n}^{r-1}(G,K),S_{n}^{r}(G,K))$ as a strong NDR pair. 
For this to be true conditions $(1)$-$(4)$ of the definition 
of an NDR pair and the extra condition of a strong NDR 
pair need to be verified. Indeed,
\begin{enumerate}
\item Note that $\tilde{u}_{r,n}^{-1}(0)=S_{n}^{r-1}(G,K)
\cap u_{r,n}^{-1}(0)$, but $u_{r,n}^{-1}(0)=F_{r}G^{n}$ and 
therefore $\tilde{u}_{r,n}^{-1}(0)=S_{n}^{r-1}(G,K)\cap 
F_{r}G^{n}=S_{n}^{r}(G,K)$.
\item Since $h_{n}(x,0)=x$ for all $x\in G^{n}$, then trivially
 $\tilde{h}_{n}(x,0)=x$ for all $x\in S_{n}^{r-1}(G,K)$.
\item $h_{n}(a,t)=a$ for all $a\in F_{r}G^{n}$, since 
$S_{n}^{r}(G,K)\subset F_{r}G^{n}$ it follows  that 
$\tilde{h}_{n}(a,t)=a$ for all $a\in S_{n}^{r}(G,K)$.
\item If $w\in \tilde{u}_{r,n}^{-1}(\left[0,1\right))=
u_{r,n}^{-1}(\left[0,1\right))\cap S_{n}^{r-1}(G,K)$, 
then $w\in u_{r,n}^{-1}(\left[0,1\right))=F_{r}G^{n}$ 
and thus $w\in F_{r}G^{n}\cap S_{n}^{r-1}(G,K)=S_{n}^{r}(G,K)$ 
and the converse is also true.  Thus, 
$\tilde{u}_{r,n}^{-1}(\left[0,1\right))=S_{n}^{r}(G,K)$.
\end{enumerate}
From the above $(\tilde{h}_{n},\tilde{u}_{r,n})$ is an NDR 
representation of $(S_{n}^{r-1}(G,K),S_{n}^{r}(G,K))$. 
Suppose now that $(x_{1},...,x_{n})\in S_{n}^{r-1}(G,K)$ is such 
that $u(x_{1},...,x_{n})<1$.
By definition 
\[
u_{r,n}(x_{1},...,x_{n})=\frac{1}{r}
\min_{1\le i_{1}<\cdots< i_{r}\le n}\left(u(x_{i_{1}})+
\cdots+u(x_{i_{r}})\right).
\]
It follows that for some $r$-tuple $1\le j_{1}<\cdots< j_{r}\le n$
\[
\frac{1}{r}\left(u(x_{j_{1}})+\cdots+u(x_{j_{r}})\right)<1.
\]
Since $0\le u(x)\le 1$ for all $x\in X$, this happens if 
and only if $u(x_{j_{s}})<1$ for some $1\le s\le r$. 
Since $(h,u)$ is a strong NDR pair representation of 
$(G,\{1_{G}\})$ then $u(h(x_{j_{s}},t))<1$ for all $t\in [0,1]$. 
On the other hand, by definition 
\[
u(h((x_{1},...,x_{n}),t))=\frac{1}{r}\min_{1\le i_{1}<
\cdots<i_{r}\le n}
\left(u(h(x_{i_{1}},t_{i_{1}}))+\cdots+
u(h(x_{i_{r}},t_{i_{r}}))\right).
\]
Here $t_{1},...,t_{n}$ are as explained before. 
For the $r$-tuple $1\le j_{1}<\cdots<j_{r}\le n$ 
\[
\frac{1}{r}(u(h(x_{j_{1}},t_{j_{1}}))+\cdots+
u(h(x_{j_{r}},t_{j_{r}})))<1
\]
as $u(h(x_{j_{s}},t_{j_{s}}))<1$. Thus in particular
\[
u(h((x_{1},...,x_{n}),t))<1.
\] 
To finish, note that $(S_{n}^{r-1}(G,K),S_{n}^{r}(G,K))$ is 
$G$-pair and if $(h,u)$ is a $G$-NDR pair representation of 
$(G,\{1_{G}\})$, then as noted above each $h_{n}$ is 
$G$-equivariant and $u_{n,r}$ is $G$-invariant. This shows that 
$(\tilde{h}_{n},\tilde{u}_{r,n})$ is in this case a $G$-NDR 
representation of $(S_{n}^{r-1}(G,K),S_{n}^{r}(G,K)).$
\qed
\medskip

Next it is shown that if $G$ is any compact Lie group and 
$K\subset Z(G)$ is a closed subgroup then such a $G$-equivariant 
NDR representation $(h,u)$ satisfying the additional properties 
required in the last proposition can always be found. 
Recall that given any Lie group $G$, if $\g$ denotes the Lie 
algebra of $G$, then the adjoint representation 
gives a homomorphism
\[
Ad:G\to Aut(\g).
\]
For each $x\in G$, $Ad_{x}$ is the derivative of the the map 
\begin{align*}
\tau_{x}:G&\to G\\
g&\mapsto xgx^{-1}.
\end{align*}
The Lie group $G$ can be seen as a Riemannian manifold. In 
particular $\g$ is a normed space with norm $\left\|\cdot \right\|$. 
In this way $Ad$ can be seen as a continuous family 
$\{Ad_{x}\}_{x\in G}$ of bounded operators 
on the normed space $\g$ that is parametrized by $G$. Since $G$ is 
assumed to be a compact Lie group, there exists 
a Haar measure $\mu$ on $G$. As usual, the 
average of the norm $\left\| \cdot \right\|$ can be used to 
obtain a new norm $\left\|\cdot \right\|_{\g} $ in $\g$ with the 
further property that each $Ad_{x}$ is an isometry 
for all $x\in G$. Indeed, define
\[
\left\|v\right\|_{\g}=\frac{1}{\mu(G)}\int_{x\in G}
\left\|Ad_{x}v \right\| d\mu(x),
\]
then $\left\|\cdot\right\|_{\g}$ is a norm on $\g$ 
that satisfies the required additional property.
From now on such a norm on $\g$ is fixed. Note that in the 
particular case when $G$ is semisimple and compact, 
the norm $\left\|\cdot \right\|_{\g}$ can be taken to be 
the negative of the Killing form which is positive definite 
and non-degenerate. 

On the other hand, for any Lie group the exponential map is a 
local homeomorphism, in particular there exists an $\epsilon >0$ 
such that the restriction of the exponential map to 
$\bar{B}_{\epsilon}(0)$
\[
\exp:\bar{B}_{\epsilon}(0)\to G
\]
is a homeomorphism onto its image. Moreover, since $K\subset G$ 
is a closed subgroup, $\epsilon$ can be chosen small enough 
so that $\exp(B_{\epsilon}(0)\cap \k)=
\exp(B_{\epsilon}(0))\cap K$, where $\k$ is the Lie 
algebra of $K$. This is used to prove the following lemma.

\begin{lemma}\label{exponential}
Let $G$ be a compact Lie group and $K\subset Z(G)$ a closed 
subgroup. Choose $\epsilon>0$ as above and let $y\in 
\bar{B}_{\epsilon/2}(0)-\{0\}$. Then for any $0<t\le 1$  
\begin{enumerate}
\item $Z_{G}(\exp(y))=Z_{G}(\exp(ty))$ and
\item $\{x\in G ~|~ [x,\exp(y)]\in K\}
=\{x\in G ~|~ [x,\exp(ty)]\in K\}$.
\end{enumerate}
\end{lemma}
\Proof
Part $(1)$ was proved in \cite[Lemma 8.1]{AC}. To prove $(2)$, 
suppose that  $c\in K\subset Z(G)$. Take $x\in G$ with 
$[x,\exp(y)]=c$; that is, $x(\exp(y))x^{-1}(\exp(-y))=c$. Since 
\[
x(\exp(y))x^{-1}=\exp(Ad_{x}y),
\] 
it follows that $\exp(Ad_{x}y)\exp(-y)=c\in Z(G)$. In particular 
$\exp(Ad_{x}y)$ and $\exp(-y)$ commute. Indeed, if $z,w\in G$ are 
such that $zw=c$ for $c$ central, then $w=z^{-1}c=cz^{-1}$ which 
implies that $zw=c=wz$. By part $(1)$ it follows that 
$\exp(Ad_{x}ty)$ and $\exp(-ty)$ commute for all $0\le t\le 1$. 
By \cite[Lemma II 2.1]{Helgason} 
\[
1=\exp(ty)\exp(-Ad_{x}ty)\exp(-ty)\exp(Ad_{x}ty)=
\exp(-t^{2}[y,Ad_{x}y]+O(t^{3})),
\]
where $(1/t^{3})O(t^{3})$ is bounded and analytic for $t$ 
sufficiently small. This shows that $[y,Ad_{x}y]=0$. Therefore
\[
c=\exp(Ad_{x}y)\exp(-y)=\exp(Ad_{x}y-y).
\]
Note that $(Ad_{x}y-y)\in \g$ is such that 
\[
\left\|Ad_{x}y-y\right\|_{\g}\le \left\|Ad_{x}-I\right\|_{\g}
\cdot \left\|y\right\|_{\g}\le (\left\|Ad_{x}\right\|_{\g}+
\left\|I\right\|_{\g})
\frac{\epsilon}{2}\le \epsilon.
\]  
The last inequality follows from the fact that 
$\left\|Ad_{x}\right\|_{\g}=1$ as $Ad_{x}$ is an isometry with 
respect to the norm $\left\|\cdot\right\|_{\g}$. 
This  means that $(Ad_{x}y-y)\in \bar{B}_{\epsilon}(0)$. 
Since $c\in K$ and $c=\exp(Ad_{x}y-y)$
with $(Ad_{x}y-y)\in \bar{B}_{\epsilon}(0)$, then as $\exp$ is 
injective on $\bar{B}_{\epsilon}(0)$ it follows that 
$k=(Ad_{x}y-y)\in \k$ where $\k$ is the Lie algebra of $K$. 
In particular, 
\[
\exp(t(Ad_{x}y-y))=\exp(tk)\in K \text{ for } t\in \R. 
\]
Unraveling 
the definitions, this means that
\[
[x,\exp(ty)]\in K \text{ for all $t\in \R$}.
\]  
Since $c\in K$ was arbitrary this shows that 
\[
\{x\in G ~|~ [x,\exp(y)]\in K\}\subset\{x\in G ~|~ 
[x,\exp(ty)]\in K\}.
\]
If $t\ne 0$, then the same argument with $s=1/t$ shows the 
other inclusion. 
Hence
\[
\{x\in G ~|~ [x,\exp(y)]\in K\}=\{x\in G ~|~ [x,\exp(ty)]\in K\}.
\]
\qed
\medskip

To finish this section, it is shown that if $G$ is any compact Lie 
group and $K\subset Z(G)$ is a closed subgroup then the previous 
lemma can be used to construct a strong $G$-equivariant NDR 
representation $(h,u)$ for the pair $(G,\{1_{G}\})$ that satisfies 
the additional properties required in Proposition 
\ref{NDR centralizers}. Take $\epsilon >0$ as before. Define a 
function
$u:G\to [0,1]$
as follows
\begin{equation*}
u(g)=\left\{ 
\begin{array}{ll}
2\left\|y\right\|_{\g}/\epsilon & \text{if }g
=\exp (y)\text{ for }g\in \exp (\bar{B}_{\epsilon /2}(0)), \\ 
1 & \text{if }g\in G-\exp (B_{\epsilon /2}(0)).%
\end{array}%
\right. 
\end{equation*}
Note that $u$ is continuous with $u^{-1}(0)=\{1_{G}\}$ and 
$u^{-1}([0,1))=\exp(\bar{B}_{\epsilon/2}(0))$ and that $u$ is 
invariant under the conjugation action of $G$. To see this, 
suppose $x\in G$, then since $Ad_{x}$ is an isometry for the 
norm $\left\|\cdot\right\|_{\g}$ 
\begin{equation*}
u(xgx^{-1})=\left\{ 
\begin{array}{ll}
2\left\|Ad_{x}(y)\right\|_{\g}/\epsilon =
2\left\|y\right\|_{\g}/\epsilon& \text{if }g=\exp (y)
\text{ for }g\in \exp (B_{\epsilon /2}(0)), \\ 
1 & \text{if }g\in G-\exp (B_{\epsilon /2}(0)).
\end{array}%
\right. 
\end{equation*}
On the other hand, 
take $s:G\to [0,1] $ any $G$-invariant continuous function 
satisfying the following properties 
\begin{equation*}
s(g)=\left\{ 
\begin{array}{ll}
1 & \text{if }g=\exp (y)\text{ for }g\in 
\exp (\bar{B}_{\epsilon /2}(0)), \\ 
0 & \text{if }g\in G-\exp (B_{\epsilon}(0)).
\end{array}%
\right. 
\end{equation*}
Clearly such a bump function $s$ always exists and can be 
constructed using partitions of unity or by hand using the 
exponential map as done in \cite{AC}. The map $s$ can be made 
to be $G$-invariant by the usual averaging trick using a Haar 
measure as before. Finally, define a homotopy 
\[
h:G\x [0,1]\to G
\] 
by 
\begin{equation*}
h(g,t)=\left\{ 
\begin{array}{ll}
\exp((1-t)y) & \text{if }g=\exp (y)\text{ for }
y\in \bar{B}_{\epsilon /2}(0), \\ 
\exp((1-s(g)t)y) & \text{if }g=\exp (y)\text{ for }y\in 
\bar{B}_{\epsilon
}(0)-B_{\epsilon /2}(0), \text{ and}\\ 
g & \text{if }g\in G-\exp (B_{\epsilon}(0)).%
\end{array}%
\right. 
\end{equation*}
The map $h$ is $G$-equivariant. To see this, take $x\in G$, 
then as $Ad_{x}$ is an isometry with respect to the metric 
$\left\|\cdot \right\|_{\g}$ the following are true
\[
h(xgx^{-1},t)=\exp((1-t)Ad_{x}y)=x\exp((1-t)y)x^{-1}
=xh(g,t)x^{-1}
\]
if  $g=\exp (y)$ for $y\in \bar{B}_{\epsilon /2}(0)$,
\[
h(xgx^{-1},t)=\exp((1-s(xgx^{-1})t)Ad_{x}y)
=x\exp((1-s(x)t)y)x^{-1}=xh(g,t)x^{-1}
\]
if $g=\exp (y)$ for $y\in \bar{B}_{\epsilon}(0)-B_{\epsilon /2}(0)$, 
and
\[
h(xgx^{-1},t)=xgx^{-1}=xh(g,t)x^{-1}
\]
if $g\in G-\exp (B_{\epsilon}(0))$. From here it follows 
that $h$ is $G$-equivariant, where $G$ is 
acting on itself by conjugation. The previous remarks can be used 
to prove the following proposition.

\begin{proposition}\label{extra}
Let $G$ be a compact Lie group and $K\subset Z(G)$ a 
closed subgroup. Then the pair $(h,u)$ as defined above 
is a strong $G$-NDR pair representation of $(G,\{1_{G}\})$ 
that satisfies the following additional properties: 
for each $g\in G-\{1_{G}\}$ and any $0\le t<1$ 
\begin{enumerate}
\item $Z(g)=Z(h(g,t))$ and
\item $\{x\in G ~|~ [x,g]\in K\}=\{x\in G ~|~ [x,h(g,t)]\in K\}$.
\end{enumerate}
\end{proposition}
\Proof
It can be seen directly, as it was done in \cite[Proposition 8.2]{AC} 
that $(h,u)$ is an NDR pair representation for $(G,\{1_{G}\})$. 
By the previous remarks, $h$ is $G$-equivariant and $u$ is 
$G$-invariant, thus $(h,u)$ is a $G$-NDR pair representation for 
$(G,\{1_{G}\})$. Moreover, $(h,u)$ is a strong $G$-NDR pair 
representation for $(G,\{1_{G}\})$. To see this, suppose that 
$g\in G$ is such that $u(g)<1$. It follows that $g=\exp(y)$ 
for some $y\in B_{\epsilon /2}(0)$. Therefore $h(g,t)=\exp((1-t)y)$ 
and 
\[
u(h(g,t))=u(\exp((1-t)y))=2\left|1-t\right| 
\left\|y\right\|_{\g}/\epsilon\le 2\left\|y\right\|_{\g}/\epsilon<1,
\]
for all $t\in [0,1]$.  Conditions $(1)$ and $(2)$ are now verified. 
To do so, note that for any $t\in [0,1]$ and any $g\in G$, $h(g,t)$ 
is either  $1_{G}$, $g$ or $\exp(ky)$ for $0<k\le 1 $ and 
$y\in \bar{B}_{\epsilon}(0)$ such that $g=\exp(y)$. Moreover, if 
$g\ne 1$ and $t\ne 1$, then  $h(g,t)\ne 1$ and thus $h(g,t)$ is 
either $g$ or $\exp(ky)$ for $0<k\le 1$. The proposition 
follows then by Lemma \ref{exponential}. 
\qed
\medskip

The following theorem which is the goal of this section is an 
immediate corollary of Propositions \ref{NDR centralizers} and 
\ref{extra}.

\begin{theorem}\label{NDR compact}
If $G$ is a compact Lie group and $K\subset Z(G)$ is a 
closed subgroup then, for each $n\ge 0$ and each $0\le r\le n$ 
the pair $(S_{n}^{r-1}(G,K),S_{n}^{r}(G,K))$ is a strong 
$G$-NDR pair. In particular,  $(\bar{S}_{n}^{r-1}(G,K),
\bar{S}_{n}^{r}(G,K))$ and $(S_{n}^{r-1}(G,K)^{H},S_{n}^{r}(G,K)^{H})$ 
are strong NDR pairs for every subgroup $H\subset G$.
\end{theorem}

\begin{remark}
The previous theorem can also be proved as follows.
Let $G$ be a compact Lie group and $K\subset Z(G)$ a closed 
subgroup. Then using the methods in \cite{ParkSuh} it 
can be proved that the spaces of the form $S_{n}^{r}(G,K)$ 
are semi-algebraic sets and that these have the homotopy type 
of a $G$-CW complex (see \cite{ParkSuh} for definitions). 
Moreover, $S^{r}_{n}(G,K)$ can be seen as a $G$-subcomplex 
of $S^{r-1}_{n}(G,K)$ for all $1\le r\le n$. The authors 
would like to thank one of the referees for pointing out
this alternative approach. The arguments provided 
here can be extended in some situations to groups which are 
not necessarily compact (see \cite{AC} for example).
\end{remark}

\section{Stable decompositions}\label{stable}

In \cite{AC}, it was proved that if $G$ is a closed 
subgroup of $GL(n,{\mathbb{C}})$, in particular if $G$ is a compact 
Lie group, there is a natural homotopy equivalence
\[
\Sigma(\Hom(\Z^{n},G))\simeq \bigvee_{1\le r\le n}\Sigma
\left(\bigvee^{\binom{{n}}{{r}}}\Hom(\Z^{r},G)/S_{r}(G) \right).
\]
Here $S_{r}(G)\subset \Hom(\Z^{r},G)$ is the subspace of $r$-tuples 
$(x_{1},...,x_{r})\in \Hom(\Z^{r},G)$ for which at least one of 
the $x_{i}$'s equals $1_{G}$. In this section it is proved that a 
similar stable splitting holds for the spaces $B_{n}(G,K)$, 
when $G$ is a compact Lie group and  $K\subset Z(G)$ is a closed 
subgroup.  

\medskip

Let $X_{*}$ be a simplicial space with face and degeneracy maps 
$\partial_{i}: X_{n}\to X_{n-1}$ and $s_{j}:X_{n}\to X_{n+1}$ 
respectively. The different degeneracy maps can be used to obtain 
a filtration of every  space $X_{n}$ as follows. 
For each $0\le r\le n$, denote
\[
S^{r}(X_{n})=\bigcup_{0\le j_{r}<\cdots<j_{1}\le n}
s_{j_{1}}s_{j_{2}}\cdots s_{j_{r}}(X_{n-r}).
\] 
Also define $S^{n+1}(X_{n})$ to be the empty set. 
These spaces form a filtration of $X_{n}$
\[
\emptyset=S^{n+1}(X_{n})\subset S^{n}(X_{n})\subset 
\cdots\subset S^{1}(X_{n})\subset S^{0}(X_{n})=X_{n}.
\]

\begin{definition}
A simplicial space $X_{*}$ is said to be simplicially NDR 
if for every $n\ge 0$ and every $0\le r\le n$ the pair 
$(S^{r}(X_{n}),S^{r+1}(X_{n}))$ is an NDR pair.
\end{definition}

Under some circumstances the previous filtration splits stably. 
For example, the following theorem proved in \cite{ABBCG} 
establishes that if a simplicial space $X_{*} $ is simplicially 
NDR, then previous filtration of $X_{n}$ splits after one 
suspension. 

\begin{theorem}\label{Adem1}
Consider a simplicial space $X_{*}$ that is simplicially NDR. 
Then for every $n\ge 0$ there is a natural homotopy equivalence
\[
\Theta(n):\Sigma(X_{n})\to \bigvee_{0\le r\le n}
\Sigma(S^{r}(X_{n})/S^{r+1}(X_{n})).
\]
\end{theorem}

In general for a simplicial space $X_{*}$ there is a natural 
filtration 
\[
F_{0}\left|X_{*}\right|\subset F_{1}\left|X_{*}\right|\subset\cdots
\subset F_{n}\left|X_{*}\right|\subset\cdots\subset 
\left|X_{*}\right|
\]
of $\left|X_{*}\right|$, the geometric realization of $X_{*}$. 
The space $F_{j}\left|X_{*}\right|$ is precisely the image in 
$\left|X_{*}\right|$ of $\sqcup_{0\le k\le j}X^{k}\x \Delta^{k}$. 
This filtration behaves nicely if $X_{*}$ satisfies certain 
conditions. One such condition is the following definition that 
was given by May in \cite{MayGeom}.

\begin{definition}
A simplicial space $X_{*}$ is proper if each pair 
$(X_{n},S^{1}(X_{n}))$ is a strong NDR pair.
\end{definition}

When a simplicial space $X_{*}$ is simplicially NDR and proper, 
then the stable homotopy type of the factors 
$S^{r}(X_{n})/S^{r+1}(X_{n})$ that appear in the previous theorem 
can be identified in terms of the filtration $F_{n}\left|X_{*}\right|$ 
of $\left|X_{*}\right|$. More precisely, there is the following 
theorem also proved in \cite{ABBCG}.

\begin{theorem}\label{Adem2}
Let $X_{*}$ be a simplicial space that is proper and simplicially 
NDR. Then there are homotopy equivalences
\[
K(n,r):\Sigma^{n+1}(S^{r}(X_{n})/S^{r+1}(X_{n}))\to 
\bigvee_{J_{r}}\Sigma^{r+1}(F_{n-r}\left|X_{*}\right|/F_{n-r-1}
\left|X_{*}\right|).
\]
Thus by Theorem \ref{Adem1} there are natural homotopy equivalences
\[
\Theta'(n):\Sigma^{n+1}(X_{n})\to \bigvee_{0\le r\le n}
\bigvee_{J_{r}}\Sigma^{r+1}(F_{n-r}\left|X_{*}\right|/F_{n-r-1}
\left|X_{*}\right|),
\]
where $J_{r}$ runs over all possible sequences of the form 
$0\le j_{r}<\cdots <j_{1}\le n$.
\end{theorem} 

As pointed out before, the collection $B_{*}(G,K)$ forms a simplicial 
space for a general Lie group and a closed central subgroup $K$. 
By definition for $0\le r\le n$ the space $S^{r}(B_{n}(G,K))$ in 
the filtration of $B_{n}(G,K)$ given by the degeneracy maps is 
precisely $S_{n}^{r}(G,K)$. In Theorem \ref{NDR compact}, 
it was proved that if $G$ is a compact Lie group and $K\subset Z(G)$ 
is a closed subgroup then for every $n\ge 0$ and every $0\le r\le n$ 
the pair $(S_{n}^{r}(G,K),S_{n}^{r+1}(G,K))$ is a strong NDR pair. 
Thus the following is obtained as a corollary of Theorem 
\ref{NDR compact}:

\begin{corollary}\label{Simplicially NDR pair}
Let $G$ be a compact Lie group and $K\subset Z(G)$ a closed 
subgroup. Then the simplicial space $B_{*}(G,K)$ is proper and 
simplicially NDR.
\end{corollary}

A stable splitting for 
$B_{n}(G,K)$ is obtained
as a direct consequence of the previous corollary.
Moreover, the different 
quotients 
\[
S^{r}(B_{n}(G,K))/S^{r+1}(B_{n}(G,K))
\] 
can be identified in the same way as in \cite[Section 6]{AC}. To 
be more precise, for $0\le r\le n$, let $J_{n,r}$ denote the set 
of all sequences of the form 
\[
1\le m_{1}<\cdots<m_{n-r}\le n.
\]
Note that $J_{n,r}$ contains precisely $\binom{{n}}{{n-r}}
=\binom{{n}}{{r}}$ elements. Given such a sequence, there is an 
associated projection
\begin{align*}
P_{m_{1},...,m_{n-r}}:B_{n}(G,K)&\to B_{n-r}(G,K)\\
(x_{1},...,x_{n})&\mapsto (x_{m_{1}},...,x_{m_{n-r}}).
\end{align*}
These projections are $G$-equivariant, with $G$ acting by 
conjugation and can be assembled to obtain a $G$-map
\begin{align*}
\eta_{n}:B_{n}(G,K)&\to \prod_{J_{n,r}}{B_{n-r}(G,K)/
S_{n-r}^{1}(G,K)}\\
(x_{1},...,x_{n})&\mapsto \{\bar{P}_{m_{1},...,m_{n-r}}
((x_{1},...,x_{n}))\}_{(m_{1},...,m_{n-r})\in J_{n,r}}.
\end{align*}
The map $\eta_{n}$ sends the $G$-invariant space 
$S_{n}^{r}(G,K)$ onto $\bigvee_{J_{n,r}}{B_{n-r}(G,K)/
S_{n-r}^{1}(G,K)}$ and $S_{n}^{r+1}(G,K)$ is mapped onto the 
base point. Therefore $\eta_{n}$ induces a $G$-equivariant 
continuous map 
\begin{align*}
S^{r}(B_{n}(G,K))/S^{r+1}(B_{n}(G,K))&\to 
\bigvee_{J_{n,r}}{B_{n-r}(G,K)/S_{n-r}^{1}(G,K)}\\
&=\bigvee^{\binom{{n}}{{r}}}{B_{n-r}(G,K)/S_{n-r}^{1}(G,K)}
\end{align*}
and this map is easily shown to be a $G$-equivariant 
homeomorphism. This together with  Theorems \ref{Adem1} and 
\ref{Adem2} can be used to prove the following theorem.

\begin{theorem}\label{decomposition}
Suppose that $G$ is a compact Lie group and that $K\subset Z(G)$ 
is a closed subgroup. Then for each $n\ge 1$ there is a natural 
$G$-equivariant homotopy equivalence
\[
\Theta(n):\Sigma(B_{n}(G,K))\simeq \bigvee_{1\le r\le n}
\Sigma\left(\bigvee^{\binom{{n}}{{r}}}B_{r}(G,K)/S^{1}_{r}(G,K)
\right).
\]
In particular there is natural homotopy equivalence
\[
\bar{\Theta}(n):\Sigma(\bar{B}_{n}(G,K))\simeq 
\bigvee_{1\le r\le n}\Sigma\left(\bigvee^{\binom{{n}}{{r}}}
\bar{B}_{r}(G,K)/\bar{S^{1}}_{r}(G,K) \right).
\]
\end{theorem}
\Proof
By  Theorems \ref{Adem1} and \ref{Simplicially NDR pair} and the 
previous remark it follows that $\Theta(n)$ is a homotopy equivalence. 
To see that this is in fact a $G$-equivariant homotopy 
equivalence, note that if $g\in G$, then the map conjugation by $g$ 
defines a map of simplicial spaces
\[
\tau_{g}:B_{*}(G,K)\to B_{*}(G,K).
\]
By naturality it follows that $\Theta(n)$ is a 
$G$-equivariant map, with $G$ acting by conjugation. On the other 
hand, if $H\subset G$ is a subgroup then $B_{n}(G,K)^{H}$ forms a 
simplicial space. By Theorem \ref{NDR compact} and Lemma \ref{G-NDR} 
the pair $(S_{n}^{r-1}(G,K)^{H},S_{n}^{r}(G,K)^{H})$ is a strong
 NDR pair. This means that the simplicial space $B_{*}(G,K)^{H}$ 
 is simplicially NDR and Theorem \ref{Adem1} applied to 
this simplicial space provides a homotopy equivalence
\[
\Theta(n,H):\Sigma(B_{n}(G,K)^{H})\simeq \bigvee_{1\le r\le n}
\Sigma\left(\bigvee^{\binom{{n}}{{r}}}B_{r}(G,K)^{H}/
S^{1}_{r}(G,K)^{H} \right).
\]
The map $\Theta(n,H)$ agrees by naturality with the fixed point 
set map $\Theta(n)^{H}$. Therefore $\Theta(n)$ is a $G$-equivariant 
map such that for every subgroup $H\subset G$, the fixed point 
map $\Theta(n)^{H}$ is a homotopy equivalence. In particular 
$\Theta(n)$ is a $G$-equivariant weak homotopy equivalence. 
By Proposition \ref{G-CW complex structure} below the $G$-space 
$B_{n}(G,K)$ has the homotopy type of a $G$-CW complex 
(similarly each $B_{r}(G,K)/S^{1}_{r}(G,K)$ has the homotopy type 
of a $G$-CW complex  for $1\le r\le n$) and thus the theorem follows
by the equivariant Whitehead Theorem. 
\qed
\medskip

\begin{remark} The naturality of the map $\Theta(n)$ in the previous 
theorem is with respect to morphisms between pairs $(G,K)$ and 
$(G',K')$; that is, $\Theta(n)$ is natural with respect to 
homomorphisms $f:G\to G'$ of Lie groups such that $f(K)\subset K'$. 
\end{remark}

\begin{proposition}\label{G-CW complex structure}
Let $G$ be a compact Lie group and $K\subset Z(G)$ a closed subgroup. 
Then for every $n\ge 1$ the space $B_{n}(G,K)$ 
has the homotopy type of a $G$-CW complex.
\end{proposition}

\Proof
The case $K=\{1_{G}\}$ is considered first, this corresponds to the 
space of commuting $n$-tuples in $G$. Since $G$ is a compact Lie 
group it is well known that $G$ carries the structure of a real 
algebraic group. This is why $B_{n}(G,\{1_{G}\})=
\Hom(\Z^{n}\,G)\subset G^{n}$ is a compact real algebraic 
$G$-variety, where the action of $G$ is given by conjugation. 
By \cite[Theorem 1.3]{ParkSuh} it follows that $\Hom(\Z^{n}\,G)$ 
has the homotopy type of a $G$-CW complex. Suppose now that
$K\subset Z(G)$ is any closed subgroup. By the previous argument, 
the space $\Hom(\Z^{n},G/K)$ has the homotopy type of a 
$G/K$-CW complex. The natural map $p:G\to G/K$ induces a $G$ action 
on $\Hom(\Z^{n},G/K)$ and this way $\Hom(\Z^{n},G/K)$ 
has the homotopy type of a $G$-CW complex. By Lemma \ref{covering} 
the natural map $\phi_{n}:B_{n}(G,K)\to \Hom(\Z^{n},G/K)$
is a $G$-equivariant locally trivial principal $K^{n}$-bundle. 
The group $G$ acts trivially on $K^{n}$ as $K$ is central and 
this is enough to conclude that $B_{n}(G,K)$ has the homotopy type 
of a $G$-CW complex. To see this, note that $B_{n}(G,K)$ is a 
separable metric space and thus by \cite[Theorem 14.2]{Murayama} 
it is enough to show that $B_{n}(G,K)$ has the homotopy type of a 
$G$-ANR. By \cite[Theorem 9.5]{Murayama} this is equivalent to 
showing that every point $\ux\in B_{n}(G,K)$ has a 
$G_{\ux}$-neighborhood which is a $G_{\ux}$-ANR but the 
latter follows easily from the local triviality of the 
principal $K^{n}$-bundle $B_{n}(G,K)\to \Hom(\Z^{n},G/K)$ 
and the fact that $G$ acts trivially on $K^{n}$.
\qed
\medskip

Using \cite[Lemma 11.3]{MayGeom} and the fact that $B_{*}(G,K)$ 
is a proper simplicial space when $G$ is a compact Lie group,
the stable factors appearing in the previous Theorem 
\ref{decomposition} can be described in terms of the natural 
filtration of $B(G,K)$, the geometric realization of 
$B_{*}(G,K)$. 

\begin{proposition}
Let $G$ be a compact Lie group and $K\subset Z(G)$ a closed 
subgroup. Then there is a $G$-equivariant homeomorphism
\[
\Sigma^{n}B_{n}(G,K)/S^{1}_{n}(G,K)\cong 
F_{n}B(G,K)/F_{n-1}B(G,K),
\]
where 
\[
F_{0}B(G,K)\subset F_{1}B(G,K)\subset\cdots\subset 
F_{n}B(G,K)\subset\cdots \subset B(G,K)
\]
is the natural filtration of $B(G,K)$.
\end{proposition}

\section{Representation spaces and symmetric products}
\label{symmetric products}
 
In this section the representation spaces $\Rep(\Z^{n},G)$ are 
studied for compact Lie groups. In particular the stable splitting 
of Theorem \ref{decomposition} is completely determined for these 
spaces when $G$ is such that $\Rep(\Z^{n},G)$ is path-connected for 
all $n\ge 1$. Also, a connection between representation spaces and 
symmetric products is explored from which interesting consequences 
can be derived.

\medskip 

To begin, take $G$ to be a compact Lie group. If $K=\{1_{G}\}$ 
is the trivial group then by definition
$B_{n}(G,\{1_{G}\})=\Hom(\Z^{n},G)$
and Theorem \ref{decomposition} provides a stable splitting for 
$\Hom(\Z^{n},G)$. This splitting agrees with \cite[Theorem 1.6]{AC} 
with the additional fact that the splitting map $\Theta(n)$ is a 
$G$-equivariant homotopy equivalence. In particular, after passing 
to orbit spaces, $\bar{\Theta}(n)$ defines a natural homotopy 
equivalence
\[
\bar{\Theta}(n):\Sigma\Rep(\Z^{n},G)\simeq \bigvee_{1\le r\le n}
\Sigma\left(\bigvee^{\binom{{n}}{{r}}}
\Rep(\Z^{n},G)/\bar{S}^{1}_{r}(G) \right).
\]
The following theorem identifies the stable pieces 
$\Rep(\Z^{n},G)/\bar{S}^{1}_{r}(G)$ under the assumption that 
$\Rep(\Z^{r},G)$ is path-connected for all $r\ge 1$.

\begin{theorem}\label{the case of rep}
Let $G$ be a compact, connected Lie group such that 
$\Rep(\Z^{r},G)$ is path-connected for $1\le r\le n$. Let $T$ be 
a maximal torus of $G$ and $W$ the Weyl group associated to $T$. 
Then there is a homeomorphism
\[
\Rep(\Z^{n},G)\cong T^{n}/W,
\]
where $W$ acts diagonally. Moreover, the map $\bar{\Theta}(n)$ 
defines a homotopy equivalence
\[
\Sigma \Rep(\Z^{n},G) \simeq \bigvee_{1\le r\le n}\Sigma
(\bigvee^{\binom{{n}}{{r}}} T^{\wedge r}/W),
\]
where $T^{\wedge r}$ is the smash product of $r$ copies of $T$.
\end{theorem}
\Proof
Consider the continuous map 
\begin{align*}
\varphi_{n}:G\times T^{n}&\to \Hom(\Z^{n},G)\\
(g,t_{1},...,t_{n})&\mapsto (gt_{1}g^{-1},...,gt_{n}g^{-1}).
\end{align*}
Under the given hypothesis, \cite[Lemma 4.2]{Baird} shows that 
every commuting $n$-tuple in $G$ lies in a maximal torus of $G$. 
Since any two maximal tori in $G$ are conjugated then  
$\varphi_{n}$ is surjective. Moreover, $\varphi_{n}$ is 
invariant under the action of $N(T)$ and therefore it induces a 
continuous map
\[
\bar{\varphi}_{n}:G\times_{N(T)}T^{n}\to \Hom(\Z^{n},G).
\]
This map is $G$-equivariant, where $G$ acts by left multiplication 
on the $G$ factor of $G\times_{N(T)}T^{n}$ and by conjugation 
on $\Hom(\Z^{n},G)$. The induced map on the level of orbit spaces 
is a homeomorphism
\[
\psi_{n}:T^{n}/W\cong \Rep(\Z^{n},G).
\]
Under this homeomorphism, $\bar{S}^{1}_{r}(G)$ corresponds to 
the subspace of $T^{n}/W$ consisting of equivalence classes of 
the form $[t_{1},...,t_{n}]$ for which at least one $t_{i}=1$. 
Therefore $\psi_{n}$ defines a homeomorphism
\[
T^{\wedge n}/W\cong \Rep(\Z^{n},G)/\bar{S}^{1}_{r}(G).
\]
\qed
\medskip

\begin{remark} By \cite[Proposition 2.3]{AC} 
if a Lie group $G$ is such that every abelian subgroup of $G$ 
is contained in a path-connected abelian subgroup, then the 
spaces $\Hom(\Z^{n},G)$ and $\Rep(\Z^{n},G)$ are path-connected 
for every $n$. This is the case for $G=U(m)$, $SU(m)$ and $Sp(m)$.
\end{remark}
 
\begin{example} Let $G=SU(m)$ with $m\ge 1$. The space 
$\Rep(\Z^{n},SU(m))$ is path-connected by the previous remark.
In this case a maximal torus $T$ 
can be taken to be the subspace of diagonal matrices 
with entries in $\S^{1}$ and determinant one, thus 
$T\cong (\S^{1})^{m-1}$. The Weyl group $W$ is the symmetric group 
$\Sigma_{m}$ acting on $T$ by permuting the diagonal entries 
of a matrix in $T$. According to the previous theorem 
there is a homotopy equivalence  
\[
\Sigma \Rep(\Z^{n},SU(m))\simeq \bigvee_{1\le r\le n}\Sigma
\left(\bigvee^{\binom{{n}}{{r}}}(\underbrace{(\S^{1})^{m-1}\wedge
\cdots \wedge(\S^{1})^{m-1}}_{r-\text{copies}})/\Sigma_{m}\right),
\]
with $\Sigma_{m}$ acting diagonally. As an example take $m=2$. 
In this case the 
stable factors are given by 
\[
(\S^{1})^{\wedge r}/\Sigma_{2}\cong \S^{r}/\Sigma_{2},
\]
with the nontrivial element $\tau\in\Sigma_{2}$ acting on a point 
$(x_{0},...,x_{r})\in \S^{r}$ by 
\[
\tau\cdot (x_{0},...,x_{r})=(x_{0},-x_{1},...,-x_{r}).
\]
\end{example}

The definition of symmetric products is recalled next.

\begin{definition}
Let $X$ be a topological space. The $m$-th symmetric product 
of $X$ is defined to be the quotient space
\[
SP^{m}(X):=X^{m}/\Sigma_{m},
\]
where the symmetric group $\Sigma_{m}$ acts by 
permuting the different factors of $X$.
\end{definition}

When $m=0$, the space $SP^{0}(X)$ is defined to be a point. 
In general, $SP^{m}(X)$ can be thought of as the space of unordered 
$m$-tuples $[y_{1},...,y_{m}]$. Suppose that $X$ is a based space 
with base point $x_{0}$. In this situation, there is a natural map
\[
SP^{m}(X)\to SP^{m+1}(X)
\]
that sends the unordered $m$-tuple $[y_{1},...,y_{m}]$ 
to $[y_{1},...,y_{m},x_{0}]$. This induces a sequence
\[
X=SP^{1}(X)\to SP^{2}(X)\to \cdots \to SP^{m}(X)\to \cdots
\]
and $SP^{\infty}(X)$ is defined as the colimit of this sequence. 
Notice that $SP^{\infty}(X)$ is precisely the free abelian 
monoid generated by $X$.

\medskip

Symmetric products naturally appear in the study of spaces of 
representations; this is explained in the following proposition. 

\begin{proposition}\label{symmetric}
There are  homeomorphisms
\begin{align*}
\Rep(\Z^{n},U(m))&\cong SP^{m}((\S^{1})^{n}),\\
\Rep(\Z^{n},Sp(m))&\cong SP^{m}((\S^{1})^{n}/\Z/2),
\end{align*}
where $\Z/2$ acts by complex conjugation on $\S^{1}$ and 
diagonally on the torus $(\S^{1})^{n}$. These homeomorphisms 
are compatible with the standard inclusions 
\[
U(m)\to U(m+1) \text{ and } Sp(m)\to Sp(m+1).
\]
\end{proposition}

\Proof
By Theorem \ref{the case of rep} when $G$ is a compact Lie group 
such that $\Rep(\Z^{n},G)$ is path-connected for all $n\ge 1$ there 
is a homeomorphism
\[
\Rep(\Z^{n},G)\cong T^{n}/W,
\]
where $T\subset G$ is a maximal torus and $W$ is the Weyl group. As 
pointed out before $\Rep(\Z^{n},G$) is path-connected for every 
$n\ge 1$ when $G=U(m)$ or $Sp(m)$. For $G=U(m)$ a maximal torus 
$T$ is of rank $m$ and can be taken to be the space of 
diagonal matrices with entries in $\S^{1}$ and the Weyl 
group $W=\Sigma_{m}$ acts by permuting the diagonal entries. 
Therefore 
\[
\Rep(\Z^{n},U(m))\cong T^{n}/W\cong ((\S^{1})^{m})^{n}/\Sigma_{m}.
\]
Here $\Sigma_{m}$ acts diagonally on the previous product 
and permutes the diagonal factors in each factor $(\S^{1})^{m}$. 
Note that 
\[
((\S^{1})^{m})^{n}/\Sigma_{m}\cong ((\S^{1})^{n})^{m}/\Sigma_{m},
\]
where $\Sigma_{m}$ acts by permuting the $(\S^{1})^{n}$ factors in 
the product $((\S^{1})^{n})^{m}$. Therefore  
\[
\Rep(\Z^{n},U(m))\cong ((\S^{1})^{n})^{m}/\Sigma_{m}\cong 
SP^{m}((\S^{1})^{n}).
\]
On the other hand, for $G=Sp(m)$ a maximal torus $T$ is homeomorphic 
to $(\S^{1})^{m}$ and the Weyl group  $W$ is a semi-direct product 
\[
W=\Sigma_{m}\ltimes (\Z/2)^{m},
\]
where $\Sigma_{m}$ acts permuting the factors in $(\S^{1})^{m}$ and 
if $\Z/2$ acts on $\S^{1}$ by complex conjugation,  
then given $(\tau_{1},...,\tau_{m})\in (\Z/2)^{m}$ and 
$(x_{1},...,x_{m})\in (\S^{1})^{m}$ then  
\[
(\tau_{1},...,\tau_{m})\cdot (x_{1},...,x_{m})=
(\tau_{1}\cdot x_{1},...,\tau_{m}\cdot x_{m}).
\]
Therefore
\[
\Rep(\Z^{n},Sp(m))\cong T^{n}/W\cong ((\S^{1})^{m})^{n}/W,
\]
with $W$ acting diagonally. In this case
\[
((\S^{1})^{m})^{n}/W\cong ((\S^{1})^{n})^{m}/W\cong 
((\S^{1})^{n}/\Z/2)^{m}/\Sigma_{m} \cong SP^{m}((\S^{1})^{n}/\Z/2).
\]
\qed
\medskip

\begin{remark} The quotient space $(\S^{1})^{n}/\Z/2$ is a singular 
space with $2^{n}$ isolated singularities. Neighborhoods around 
the singular points look like cones on $\R P^{n-1}$. For example, 
when $n=4$ the quotient $(\S^{1})^{4}/\Z/2$ 
is an orbifold with $16$ isolated singularities.  
\end{remark}

The identification of $\Rep(\Z^{n},G)$ as a symmetric product for 
$G=U(m)$ and $Sp(m)$ has interesting consequences. For 
example the following can be derived. 

\begin{corollary}
For every $m\ge 1$ there is a homeomorphism
\[
\Rep(\Z^{2},Sp(m))\cong \C P^{m}.
\]
\end{corollary}
\Proof
By the previous corollary 
\[
\Rep(\Z^{2},Sp(m))\cong SP^{m}((\S^{1}\times \S^{1})/\Z/2).
\]
The quotient space $(\S^{1}\times \S^{1})/\Z/2$ is homeomorphic to 
$\S^{2}$. Therefore 
\[
\Rep(\Z^{2},Sp(m))\cong SP^{m}(\S^{2})\cong \C P^{m}.
\]
\qed
\medskip

Suppose now that $G$ is a topological abelian group. The 
multiplication on $G$ defines a continuous function
\begin{align*}
\mu_{m}:SP^{m}(G)&\to G\\
[x_{1},...,x_{m}]&\mapsto x_{1}\cdots x_{m}.
\end{align*}  
When $G=(\S^{1})^{n}$ the following diagram commutes
\[
\xymatrix{
SP^{m}((\S^{1})^{n})\ar[d]_{\cong}\ar[r]^{\mu_{m}}&(\S^{1})^{n}\\
\Rep(\Z^{n},U(m))\ar[ru]_{det_{*}}&
}
\]
and thus $\mu_{m}$ can be identified with the map induced by the 
determinant $det_{*}$. Note that $det_{*}^{-1}(1,...,1)$ is 
precisely the space of commuting $n$-tuples of elements in $SU(m)$ 
modulo conjugation in $U(m)$. This space agrees with 
$\Rep(\Z^{n},SU(m))$ because any two elements 
in $\Hom(\Z^{n},SU(m))$ are conjugate by an element in $U(m)$ if 
and only if they are conjugate by an element in $SU(m)$. Moreover 
the map 
\[
\mu_{m}:SP^{m}((\S^{1})^{n})\to (\S^{1})^{n}
\]
is a locally trivial fiber bundle.
This can be seen by picking an $n$-th root of the determinant 
function which can always be done locally. For the particular 
case of $n=2$, the map $\mu_{m}$ agrees with the Abel-Jacobi map 
\[
\mu_{m}:SP^{m}(\S^{1}\times \S^{1})\to \S^{1}\times \S^{1}.
\]
In general, if $\M_{g}$ is a Riemann surface of genus 
$g$ and $J(\M_{g})$ denotes its Jacobian then the Abel-Jacobi is 
a map 
\[
\mu_{m}:SP^{m}(\M_{g})\to J(\M_{g}),
\]
that makes $SP^{m}(\M_{g})$ into a fiber bundle 
over $J(\M_{g})$ with fiber type $\C P^{m-g}$  for $m\ge 2g$. 
In particular $\Rep(\Z^{2},SU(m))\cong \C P^{m-1}$ for 
$m\ge 2$. This proves the following proposition.

\begin{proposition}
For every $n$ and $m\ge 1$ composition with the determinant defines 
a locally trivial fiber bundle 
\[
det_{*}:\Rep(\Z^{n},U(m))\to (\S^{1})^{n}.
\] 
with fiber type $\Rep(\Z^{n},SU(m))$. In particular, for $m\ge 2$
\[
\Rep(\Z^{2},SU(m))\cong \C P^{m-1}.
\]
\end{proposition}

The identification $\Rep(\Z^{2},U(m))\cong 
SP^{m}(\S^{1}\times \S^{1})$ can also be 
used to prove the following corollary using the results 
from \cite{Macdonald}. 

\begin{corollary}\label{cohomology calculation}
\begin{equation*}
H^{i}(\Rep(\Z^{2},U(m));\Z)\cong \left\{ 
\begin{array}{ccc}
\Z & \text{ if } & i=0,\\ 
\Z\oplus \Z& \text{ if }   &1\le i\le 2m-1,  \\ 
\Z & \text{ if } & i=2m,\\
0 & \text{if} & i>2m.
\end{array}%
\right. 
\end{equation*}%
\end{corollary}

In general the cohomology groups of the spaces $\Rep(\Z^{n},U(m))$ 
and $\Rep(\Z^{n},Sp(m))$ can be computed using the identifications 
given in Proposition \ref{symmetric}. Note that by 
\cite[Section 22]{Steenrod1} the inclusions  
\[
\Rep(\Z^{n},U(m))\to \Rep(\Z^{n},U(m+1)) \text{ and } 
\Rep(\Z^{n},Sp(m))\to \Rep(\Z^{n},Sp(m+1))
\] 
split 
on the level of cohomology. 

\begin{definition}
For each $n\ge 1$, define $\Rep(\Z^{n},SU)$ to be the colimit of 
the sequence 
\[
\Rep(\Z^{n},SU(1))\to \Rep(\Z^{n},SU(2))\to \cdots \to 
\Rep(\Z^{n},SU(m))\to\Rep(\Z^{n},SU(m+1)) \cdots  
\]
induced by the canonical inclusions $SU(m)\to SU(m+1)$. Similarly, 
$\Rep(\Z^{n},U)$ and $\Rep(\Z^{n},Sp)$ are defined as the colimit 
of the finite stages under the canonical inclusions.
\end{definition}

The homotopy type of these spaces is established in the following 
theorem\footnote{A proof for the case of the unitary groups also 
appears in \cite[Theorem 6.6]{Ramras}.}

\begin{theorem}
For every $n\ge 1$ there are homotopy equivalences
\begin{align*}
\Rep(\Z^{n},SU)\simeq &\prod_{2\le i\le n}
K(\Z^{\binom{{n}}{{i}}},i),\\
\Rep(\Z^{n},U)\simeq &\prod_{1\le i\le n}
K(\Z^{\binom{{n}}{{i}}},i),\\
\Rep(\Z^{n},Sp)\simeq &\prod_{1\le i \le 
\left\lfloor n/2\right\rfloor}
K(\Z^{\binom{{n}}{{2i}}}\oplus (\Z/2)^{r(2i)},2i),
\end{align*}
where the $r(i)$'s are integer defined by
\begin{equation*}
r(i)=\left\{
\begin{array}{ccc}
\binom{{n}}{{0}}+\binom{{n}}{{1}}+\cdots+\binom{{n}}{{n-i-1}}& 
\text{ if } &1\le i \le n \\
0& \text{ if } &i> n.
\end{array}%
\right. 
\end{equation*}%
\end{theorem}
\Proof
By Proposition \ref{symmetric} there are  homeomorphisms
\begin{align*}
\Rep(\Z^{n},U(m))&\cong SP^{m}((\S^{1})^{n}),\\
\Rep(\Z^{n},Sp(m))&\cong SP^{m}((\S^{1})^{n}/\Z/2),
\end{align*}
which are compatible with the standard inclusions. It follows 
that there are homeomorphisms 
\[
\Rep(\Z^{n},U)\cong SP^{\infty}((\S^{1})^{n}) \text{ and } 
\Rep(\Z^{n},Sp)\cong 
SP^{\infty}((\S^{1})^{n}/\Z/2).
\]
In general in \cite{DoldThom} it is proved that if $X$ is a 
connected CW complex of finite type then 
\[
SP^{\infty}(X)\simeq \prod_{n\ge 1}K(H_{n}(X,\Z),n).
\]
This is then used to handled the cases of $U$ and $Sp$. 
(The homology groups of $(\S^{1})^{n}/\Z/2$ are computed below 
in proposition \ref{homology quotients}).  
To handle the case of $SU$ notice that a homotopy equivalence 
\[
f:\prod_{1\le i\le n}K(\Z^{\binom{{n}}{{i}}},i)
\to SP^{\infty}((\S^{1})^{n})\cong\Rep(\Z^{n},U)
\]
can be constructed as follows. For each $1\le i\le n$ pick 
\[
f_{i}:\bigvee^{\binom{{n}}{{i}}}\S^{i}\to 
SP^{\infty}((\S^{1})^{n})
\] 
inducing an isomorphism in $\pi_{i}$ for $2\le i\le n$ and 
an isomorphism on $H_{1}$ for $i=1$. The map 
\[
f_{1}:\bigvee^{n}\S^{1}\to SP^{\infty}((\S^{1})^{n})
\]
can be chosen to be the inclusion
\[
\bigvee^{n}\S^{1}\to (\S^{1})^{n}=SP^{1}((\S^{1})^{n})\to 
SP^{\infty}((\S^{1})^{n})
\]
where the $j$-th wedge factor in $\bigvee^{n}\S^{1}$ is mapped 
into the $j$-th factor of $(\S^{1})^{n}$. Since 
$SP^{\infty}((\S^{1})^{n})$ has the structure of a monoid the 
different $f_{i}'s$ can be assembled to get a homomorphism
\[
g:SP^{\infty}(\bigvee_{1\le i\le n}\bigvee^{\binom{{n}}{{i}}}\S^{i})
\to \Rep(\Z^{n},U).
\]
This map is a homotopy equivalence. Also there is a homotopy 
equivalence
\[
h:(\S^{1})^{n}\times\prod_{2\le i\le n}K(\Z^{\binom{{n}}{{i}}},i)
=\prod_{1\le i\le n}SP^{\infty}(\bigvee^{\binom{{n}}{{i}}}\S^{i})
\cong SP^{\infty}(\bigvee_{1\le i\le n}
\bigvee^{\binom{{n}}{{i}}}\S^{i}),
\]
and $f$ is defined to be the composition $f=g\circ h$. On the 
other hand, the locally trivial fibration sequences  
\[
\Rep(\Z^{n},SU(m))\to \Rep(\Z^{n},U(m))
\stackrel{det_{*}}{\rightarrow} (\S^{1})^{n}
\]
are compatible with the standard inclusions $SU(m)\to SU(m+1)$ and 
$U(m)\to U(m+1)$. Therefore the determinant gives rise to a well 
defined function 
\[
det_{*}:\Rep(\Z^{n},U)\to (\S^{1})^{n}
\]
which is a fibration with fiber type $\Rep(\Z^{n},SU)$. Moreover, by 
the way the map $f_{1}$ was chosen, it follows that there is 
a commutative diagram of fibrations 
\[
\begin{CD}
\prod_{2\le i\le n}K(\Z^{\binom{{n}}{{i}}},i)
@>f_{|}>>\Rep(\Z^{n},SU)\\
@VVV     @VVV\\   
(\S^{1})^{n}\times\prod_{2\le i\le n}K(\Z^{\binom{{n}}{{i}}},i)
@>f>> \Rep(\Z^{n},U)\\
@V{p_{1}}VV     @VV{det_{*}}V\\ 
(\S^{1})^{n} @>{=}>>  (\S^{1})^{n}\\
\end{CD}
\]
From the long exact sequences on homotopy groups associated to 
these fibrations and the five lemma it follows that 
\[
f_{|}:\prod_{2\le i\le n}K(\Z^{\binom{{n}}{{i}}},i)
\to \Rep(\Z^{n},SU)
\]
is a weak homotopy equivalence. The proposition follows by noting 
that $\Rep(\Z^{n},SU)$ has the homotopy type of a CW complex.
\qed
\medskip

\begin{proposition}\label{homology quotients}
The homology groups of $(\S^{1})^{n}/\Z/2$ are given by 
\begin{equation*}
H_{i}((\S^{1})^{n}/\Z/2;\Z)\cong \left\{ 
\begin{array}{ccl}
\Z &  \text{ if } & i=0,\\
\Z^{\binom{{n}}{{i}}}\oplus 
(\Z/2)^{\binom{{n}}{{0}}+\binom{{n}}{{1}}+\cdots+
\binom{{n}}{{n-i-1}}} & 
\text{ if } & i \text{ is even  and } 0< i\le n,\\ 
0 & \text{else}. & 
\end{array}%
\right. 
\end{equation*}
\end{proposition}
\Proof
Let $G=\Z/2$ and $X=(\S^{1})^{n}$ seen as a $G$-space with $G$ 
acting diagonally and by complex conjugation on each $\S^{1}$ 
factor. Denote by $F$ the fixed point set under this $G$ action 
which in this case is a 
discrete set with $2^{n}$ points. The cohomology groups 
$H^{i}(X/G;\Z)$ are computed first. Using the long exact sequence 
in cohomology associated to the pair $(X/G,F)$, it suffices to 
compute the relative cohomology groups $H^{*}(X/G,F;\Z)$. To do so, 
consider the map 
$\varphi: X\times_{G}EG\to X/G$
induced by the $G$-equivariant projection $\pi_{1}:X\times EG\to X$.
By \cite[Proposition 1.1 VII]{Bredon} the map $\varphi$ induces 
an isomorphism
\begin{equation}\label{IsomorphismBredon}
\varphi^{*}:H^{*}(X/G,F;\Z)\to H_{G}^{*}(X,F;\Z).
\end{equation}
The equivariant cohomology groups $H_{G}^{*}(X,F;\Z)$ are 
computed next. 
This is achieved by studying the map 
$j_{G}^{*}:H_{G}^{*}(X;\Z)\to H_{G}^{*}(F;\Z)$
induced by the inclusion $j:F\hookrightarrow X$.
Consider the fibration sequence 
\[
X\to X\times_{G}EG\to BG.
\]
By \cite[Theorem 1.2]{AGPP} the Lyndon-Hochshild-Serre spectral 
sequence associated to this fibration collapses on the 
$E_{2}$-term and there are no extension problems. Therefore
\[
H^{k}_{G}(X;\Z)=H^{k}(X\times_{G}EG;\Z)\cong 
\bigoplus_{i+j=k}H^{i}(G;H^{j}(X;\Z)).
\]
Using \cite[Proposition 1.10]{Adem} these cohomology groups can 
be computed explicitly to obtain
\begin{equation}\label{equivariantcohomologyX}
H^{k}_{G}(X;\Z)\cong \left\{ 
\begin{array}{ccl}
\Z &  \text{ if } & k=0,\\
\Z^{\binom{{n}}{{k}}}\oplus (\Z/2)^{\binom{{n}}{{0}}+
\cdots+\binom{{n}}{{k-1}}} & 
\text{ if } & k \text{ is even  and } 0< k\le n,\\ 
(\Z/2)^{2^{n}} & \text{ if } & k \text{ is even }, k>n,\\ 
0 & \text{ if } & k \text { is odd}. 
\end{array}%
\right. 
\end{equation}
In particular, the map $j_{G}^{*}:H_{G}^{k}(X;\Z)
\to H_{G}^{k}(F;\Z)$ is trivial for $k$ odd. On the other hand, 
since $G$ acts trivially on $F$,
then $F\times_{G}EG=F\times BG$ and thus
\begin{equation}\label{equivariantcohomologyF}
H^{k}_{G}(F;\Z)=H^{k}(F\times BG)\cong \left\{ 
\begin{array}{ccl}
\Z^{2^{n}} &  \text{ if } & k=0,\\
(\Z/2)^{2^{n}}&  \text{ if } & k \text{ is even}, k>0,\\
0& \text{ if } & k \text{ is odd}.
\end{array}%
\right. 
\end{equation}
By the universal coefficient theorem the natural map
\[
H^{k}(X\times_{G}EG;\Z)\otimes \Z/2\to H^{k}(X\times_{G}EG;\Z/2)
\]
is an isomorphism when $k$ is even since in this case the Tor term 
vanishes as
\[
H^{k+1}(X\times_{G}EG;\Z)=0 \text{ when $k$ is even}.
\] 
The same is true when $X$ is replaced by $F$. Therefore
when $k$ is even there is a commutative diagram
\begin{equation}\label{diagramUCT}
\xymatrix{
H^{k}(X\times_{G}EG;\Z)\otimes \Z/2\ar[r]^{\cong}
\ar[d]_{j_{G}^{*}\otimes\Z/2}   
&H^{k}(X\times_{G}EG;\Z/2)\ar[d]^{j_{G,\Z/2}^{*}}\\
H^{k}(F\times_{G}EG;\Z)\otimes \Z/2\ar[r]^{\cong} 
&H^{k}(F\times_{G}EG;\Z/2).
}
\end{equation}
On the other hand, since $H^{k}(X;\Z/2)=0$ for $k>n$, then 
by \cite[Theorem 1.5 VII]{Bredon} 
the map 
\[
j_{G,\Z/2}^{*}:H^{k}_{G}(X;\Z/2)\to H^{k}_{G}(F;\Z/2)
\] 
is an isomorphism for $k>n$. Moreover, for any 
$k\ge 0$ the map $j_{G,\Z/2}^{*}$ is injective. To see this note 
that 
\[
\sum_{i\ge 0}rk H^{k}(X,\Z/2)=\sum_{i\ge 0}rk H^{k}(F,\Z/2)=2^{n}.
\]
Therefore, by \cite[Theorem 1.6 VII]{Bredon} $X$ is totally 
nonhomologous to zero (mod $p$)
\footnote[1]{A $G$-space $X $ is said to be totally nonhomologous 
to zero (mod $p$) if the restriction to a typical fiber 
$H^{*}(X\times_{G} EG;\Z/p)\to H^{*}(X;\Z/p)$ is
surjective.}  in $X\times_{G}EG$ and by \cite[Theorem 1.5 VII]{Bredon} 
the map $j_{G,\Z/2}^{*}$ is injective. Notice that 
when $k$ is even by (\ref{equivariantcohomologyF}) it follows that
$H^{k}_{G}(F;\Z)\otimes \Z/2\cong H^{k}_{G}(F;\Z)$.
Using this and (\ref{diagramUCT}) it follows 
that when $k$ is even the image of 
\[
j_{G}^{*}:H_{G}^{k}(X;\Z)\to H_{G}^{k}(F;\Z)
\]
is a $\Z/2$-vector space of rank 
$\binom{{n}}{{0}}+\binom{{n}}{{1}}+\cdots+\binom{{n}}{{k}}$. 
Similarly, using (\ref{equivariantcohomologyX}), 
(\ref{equivariantcohomologyF}) and (\ref{diagramUCT}) it can be 
seen that when $k$ is even 
$\text{Ker}(j_{G}^{*})\cong \Z^{\binom{{n}}{{k}}}$. By the long exact
sequence in $G$-equivariant cohomology associated to the pair 
$(X,F)$ if follows that
\begin{equation*}
H^{k}_{G}(X,F;\Z)\cong \left\{ 
\begin{array}{ccl}
0& \text{ if }& k=0,\\
\Z^{2^{n}-1}& \text{ if }& k=1,\\
\Z^{\binom{{n}}{{k}}} &  \text{ if } & k \text{ is even}, 0<k\le n,\\
(\Z/2)^{\binom{{n}}{{0}}+\binom{{n}}{{1}}+\cdots+\binom{{n}}{{n-k}}}&
  \text{ if } & k \text{ is odd}, 1<k\le n,\\
0& \text{ if } & k>n.
\end{array}%
\right. 
\end{equation*}
Since $F$ is discrete, from the long exact sequence in 
cohomology associated to the pair $(X/G,F)$ it follows that 
$H^{k}(X/G,F;\Z)\cong H^{k}(X/G;\Z)$ for $k\ge 2$. It is easy to 
see that $H^{0}(X/G;\Z)=\Z$ and $H^{1}(X/G;\Z)=0$. Using this and 
the isomorphism (\ref{IsomorphismBredon}) the cohomology groups 
$H^{k}(X/G;\Z)$ are computed. Finally, the result follows by applying 
the universal coefficient theorem to obtain the 
homology groups $H_{i}(X/G;\Z)$.
\qed
\medskip

On the other hand, as noted before a stable splitting 
for the spaces of the from $\Rep(\Z^{n},G)$ is obtained 
by Theorem \ref{decomposition}. Via the identifications 
\[
\Rep(\Z^{n},U(m))\cong SP^{m}((\S^{1})^{n}) \text{ and } 
\Rep(\Z^{n},Sp(m))\cong SP^{m}((\S^{1})^{n}/\Z/2)
\]
a stable splitting for the spaces $SP^{m}((\S^{1})^{n})$ and 
$SP^{m}((\S^{1})^{n}/\Z/2)$ is then obtained. It turns out that
a similar splitting holds in general for spaces of the form 
$SP^{m}(G^{n})$ when $G$ is a topological group 
with nondegenerate unit. This can be seen using the following 
definition given by May \cite[Definition A.7]{MayGeom}.

\begin{definition} Let $\mathcal{U}$ be the category of compactly 
generated Hausdorff spaces. A functor $F:\mathcal{U}\to \mathcal{U}$, 
is said to be admissible if for any NDR representation of $(h,u)$ 
of $(X,A)$ as an NDR pair determines a representation $(Fh,Fu)$ 
of $(FX,FA)$ as an NDR pair such that $(Fh)_{t}=F(h_{t})$ and 
such that, for any map $g:X\to X$ that satisfies the condition 
$u(g(x))<1$ whenever $u(x)<1$, then the map $Fu:FX\to [0,1]$ 
satisfies $(Fu)(Fg)(y)<1$ whenever $Fu(y)<1$, $y\in FX$. 
\end{definition}

For every $m\ge 1$ the functor $SP^{m}$ is an admissible functor. 
Indeed, if $(X,A)$ is an NDR pair represented by $(h,u)$, then 
$(SP^{m}X,SP^{m}A)$ is an NDR pair represented by 
$(SP^{m}h,SP^{m}u)$, where 
\begin{align*}
SP^{m}h(([x_{1},..,x_{m}],t))&=[h(x_{1},t),...,h(x_{m},t)]\\
SP^{m}u([x_{1},..,x_{m}])&
=\displaystyle\mathop{\text{max}}_{1\le i\le m}u(x_{i}).
\end{align*}

The other condition of an admissible functor is easily verified. 
From the definition it follows that if $F$ is an admissible 
functor and  $X_{*}$ is a simplicial space that is simplicially 
NDR (resp. proper), then the simplicial space $FX_{*}$ is also 
simplicially NDR (resp. proper). In particular, a stable 
decomposition as in Theorem \ref{Adem1} for $FX_{n}$ is obtained 
whenever this theorem applies for $X_{n}$. Therefore if 
$G$ is a topological group with nondegenerate unit; that is; 
$(G,\{1_{G}\})$ is an NDR pair, then the simplicial space $B_{*}G$ 
is simplicially NDR and it is proper if $(G,\{1_{G}\})$ is further 
assumed to be a strong NDR pair. Let $S^{1}(SP^{m}(G^{n}))$  
be the subspace of $SP^{m}(G^{n})$ consisting of unordered 
$m$-tuples $\ux=[x_{1},...,x_{m}]$ with $x_{i}\in G^{n}$ 
with the following property. If 
\[
x_{i}=(g_{i1},...,g_{in}),
\] 
then there is a $j$ such that $g_{ij}=1_{G}$ for all 
$1\le i\le n$. If the elements $g_{ij}$ are seen as the 
entries of an $m\times n$ matrix $A(\ux)$, then this condition 
means that at least one of the columns of $A(\ux)$ has entries 
all equal to $1_{G}$. From the previous remark the 
following corollary is obtained.

\begin{corollary}
For any topological group $G$ with nondegenerate unit 
and any $n$ and $m\ge 1$ there is a natural homotopy equivalence
\[
\Sigma SP^{m}(G^{n})\simeq \bigvee_{1\le r\le n} 
\Sigma\left(\bigvee^{\binom{{n}}{{r}}} 
SP^{m}(G^{r})/S^{1}(SP^{m}(G^{r}))\right).
\]
If in addition $(G,\{1\})$ is a strong NDR pair, then each 
factor $SP^{m}(G^{r})/SP^{m}(F^{1}G^{r})$ is stably equivalent 
to the quotient $F_{r}SP^{m}BG/F_{r-1}SP^{m}BG$, for a 
filtration of $SP^{m}BG$
\[
F_{0}SP^{m}BG\subset F_{1}SP^{m}BG\subset \cdots 
\subset F_{r}SP^{m}BG\subset \cdots \subset SP^{m}BG.
\]
\end{corollary}
\Proof
From the previous remark, the simplicial space $SP^{m}(G^{*})$ 
is simplicially NDR and proper if $(G,\{1\})$ is a 
strong NDR pair. By Theorem \ref{Adem1}, the
filtration 
\[
\emptyset=S^{n+1}(SP^{m}(G^{n}))\subset S^{n}(SP^{m}(G^{n}))
\subset\cdots \subset S^{1}(SP^{m}(G^{n}))\subset SP^{m}(G^{n})
\]
of $SP^{m}(G^{n})$ given by the degeneracy maps 
splits after one suspension and there 
is a natural homotopy equivalence
\[
\Sigma SP^{m}(G^{n})\simeq \bigvee_{0\le r\le n} 
\Sigma (S^{r}(SP^{m}(G^{n}))/S^{r+1}(SP^{m}(G^{n}))).
\]
The first part of the corollary follows by showing that 
\[
S^{r}(SP^{m}(G^{n}))/S^{r+1}(SP^{m}(G^{n}))\cong 
\bigvee^{\binom{{n}}{{n-r}}} SP^{m}(G^{n-r})/S^{1}(SP^{m}(G^{n-r})).
\]
This is achieved in the same way as in the case of the 
spaces $B_{n}(G,K)$ using the different projection maps
\begin{align*}
P_{m_{1},...,m_{n-r}}:G^{n}&\to G^{n-r}\\
(x_{1},...,x_{n})&\mapsto (x_{m_{1}},...,x_{m_{n-r}}).
\end{align*}
for all possible sequences $1\le m_{1}\le \cdots m_{n-r}\le n$.
To show the last part note that if $(G,\{1_{G}\})$ is a strong 
NDR pair, then $B_{*}G$ is a proper simplicial space. By 
\cite[Lemma 11.3]{MayGeom} each factor 
$SP^{m}(G^{r})/S^{1}(SP^{m}(G^{r}))$ is stably equivalent to 
the quotient $F_{r}|SP^{m}B_{*}G|/F_{r-1}|SP^{m}B_{*}G|$, where 
\[
F_{0}|SP^{m}B_{*}G|\subset F_{1}|SP^{m}B_{*}G|\subset
\cdots \subset F_{r}|SP^{m}B_{*}G|\subset \cdots |SP^{m}B_{*}G| 
\]
is the natural filtration of $|SP^{m}B_{*}G|$. 
But in the category of compactly generated weak Hausdorff spaces 
$|SP^{m}B_{*}G|\cong SP^{m}BG$.
\qed
\medskip

\section{Compact Lie groups of rank one}

In this section the stable factors appearing in 
Theorem \ref{decomposition} are studied for the particular case 
of compact connected rank one Lie groups. 

\medskip

Let $G$ be a compact, connected rank one Lie group. 
If $G$ is such a Lie group then it is isomorphic 
to $\S^{1}$, $SU(2)$ or $SO(3)$. A similar argument 
can be used to understand the stable factors 
$B_{n}(G,K)/S_{n}^{1}(G,K)$ in these
three cases. The case $K=\{1_{G}\}$ is handled first, this 
corresponds to the space of commuting $n$-tuples in $G$. Fix $T$ 
a maximal torus for $G$ and let $W$ denote the corresponding 
Weyl group. When $G=\S^{1}$ then $T=\S^{1}$, for $G=SU(2)$ the 
maximal torus $T$ can be taken to be the space of diagonal matrices 
with entries in $\S^{1}$ and determinant one and when $G=SO(3)$, 
$T$ can be taken to be the space of rotations with respect
to the $z$-axis.  The possibilities for $W$ are the trivial group 
when $G=\S^{1}$ or $\Sigma_{2}$ acting on the Lie algebra 
$\t\cong i\R$ by complex conjugation in the other cases. In any of 
these cases there is a $W$-equivariant homeomorphism 
\[
\psi: \t\to T-\{1\}
\]
sending $0$ to $-1\in T$. When $G=\S^{1}$, $\psi$ can be taken 
to be the transformation
\begin{align*}
\psi:\t\cong i\R&\to T-\{1\}\\
z&\mapsto \frac{z+1}{z-1}.
\end{align*}
When $G=SU(2)$, the map $\psi$ can be taken to be
\begin{align*}
\psi:\t\cong i\R&\to T-\{1\}\\
z&\mapsto \left(\begin{array}{cc}
\frac{z+1}{z-1} &  \\ 
& \frac{z-1}{z+1}
\end{array}\right),
\end{align*}
and similarly for $SO(3)$.

\medskip

The main tool for understanding the stable factors 
$\Hom(\Z^{n},G)/S^{1}_{n}(G)$ is the map already used in 
Theorem \ref{the case of rep}
\begin{align*}
\bar{\varphi}_{n}:G/T\times_{W} T^{n}&\to \Hom(\Z^{n},G)_{(1,...,1)}\\
(\bar{g},t_{1},...,t_{n})&\mapsto (gt_{1}g^{-1},...,gt_{n}g^{-1}),
\end{align*}
where $\Hom(\Z^{n},G)_{(1,...,1)}$ denotes the path-connected 
component of $\Hom(\Z^{n},G)$ containing the $n$-tuple $(1,...,1)$. 
Note that trivially $\Hom(\Z^{1},G)/S_{1}^{1}(G)=G$ and thus from now 
on $n$ is assume to be at least $2$. As pointed out before, 
$\bar{\varphi}_{n}$ is a continuous surjective map. Also, if $g\in G$
and $t\in T$, then $gtg^{-1}=1$ if and only if $t=1$. This shows that  
the restriction of $\bar{\varphi}_{n}$ defines a surjective 
continuous map
\[
\alpha_{n}=(\bar{\varphi}_{n})_{|}:G/T\times_{W} (T-\{1\})^{n}
\to (\Hom(\Z^{n},G)_{(1,...,1)}-S^{1}_{n}(G)_{(1,...,1)}),
\]
where $S^{1}_{n}(G)_{(1,..,1)}=S^{1}_{n}(G)\cap 
\Hom(\Z^{n},G)_{(1,...,1)}$. 
Since the inclusion map 
\[
S^{1}_{n}(G)\hookrightarrow \Hom(\Z^{n},G)
\] 
is a cofibration, in particular
\[
S^{1}_{n}(G)_{(1,...,1)}\hookrightarrow \Hom(\Z^{n},G)_{(1,...,1)}
\] 
is also a cofibration, it follows that 
\[
\Hom(\Z^{n},G)_{(1,...,1)}/S^{1}_{n}(G)_{(1,...,1)}\cong 
(\Hom(\Z^{n},G)_{(1,...,1)}-S^{1}_{n}(G)_{(1,...,1)})^{+},
\]
where in general for a space $X$, $X^{+}$ denotes its one point 
compactification. On the other hand, using the $W$-equivariant 
homeomorphism $\psi$, 
it follows that 
\[
G/T\times_{W} (T-\{1\})^{n}\cong G/T\times_{W} \t^{n}.
\]
Note that $W$ acts freely on $G/T$ (this is the induced action of 
$G$ acting on itself on the right) and thus the projection map 
$p_{n}:G/T\times_{W} \t^{n}\to (G/T)/W$ is a vector bundle 
over $(G/T)/W=G/N(T)$. The map $\bar{\varphi}_{n}$ and its 
restriction $\alpha_{n}$ are not injective for a general Lie 
group. The different possibilities are considered next.

\medskip

When $G=\S^{1}$, the map $\alpha_{n}$ is a homeomorphism 
and $G/N(T)=*$. Since $\Hom(\Z^{n},\S^{1})$ is connected, 
then $\Hom(\Z^{n},\S^{1})/S^{1}_{n}(\S^{1})$ is 
homeomorphic to the one point compactification of $\t^{n}$; 
that is, in this case the stable factors are given by
\[
\Hom(\Z^{n},\S^{1})/S^{1}_{n}(\S^{1})\cong \S^{n}.
\] 
Of course, this can be seen trivially using the fact that 
$\Hom(\Z^{n},\S^{1})=(\S^{1})^{n}$ as $\S^{1}$ is abelian. 
Modulo conjugation the situation is 
trivial since $\S^{1}$ is abelian, thus 
\[
\Rep(\Z^{n},\S^{1})/\bar{S}^{1}_{n}(\S^{1})\cong \S^{n}.
\] 

\medskip

When $G$ equals $SU(2)$ the map $\alpha_{n}$ is no longer 
injective. Let 
\[
s_{n}:G/N(T)\to G/T\times_{W} \t^{n}
\] 
be the zero section of the vector bundle $p_{n}$. Then in this case 
\[
\alpha_{n}^{-1}(-1,...,-1)=s_{n}(G/N(T))
\]
and $\alpha_{n}^{-1}(g_{1},...,g_{n})$ is a point for 
$(g_{1},...,g_{n})\ne (-1,...,-1)$. On the other hand, the space 
$\Hom(\Z^{n},SU(2))$ is connected, $G/T\cong \S^{2}$ with $W$ acting 
as the antipodal map and therefore $G/N(T)\cong \R P^{2}$. 
If $\lambda_{n}$ denotes the canonical vector bundle over $\R P^{n}$, 
then the vector bundle $p_{n}:G/T\times_{W} \t^{n}\to (G/T)/W$ is 
precisely $n\lambda_{2}$, the Whitney sum of $n$ copies of 
$\lambda_{2}$. 
Therefore
\begin{equation*}
\Hom(\Z^{n},SU(2))/S^{1}_{n}(SU(2))\cong\left\{ 
\begin{array}{cl}
\S^{3} & \text{if }n=1, \\  
(\R P^{2})^{n\lambda_{2}}/s_{n}(\R P^{2})& \text{if }n\ge 2 
\end{array}%
\right. 
\end{equation*}
where $X^{\mu}$ denotes the associated Thom space for a vector 
bundle $\mu$ over $X$. This agrees with the computations given 
independently in  \cite{Baird2} and \cite{Crabb}. Modulo 
conjugation the situation simplifies. As pointed out before 

\[
\Rep(\Z^{n},SU(2))/\bar{S}^{1}_{n}(SU(2))\cong \S^{n}/\Sigma_{2}
\]
with the nontrivial element $\tau\in\Sigma_{2}$ acting on a point 
$(x_{0},...,x_{n})\in \S^{r}$ by 
\[
\tau\cdot (x_{0},...,x_{r})=(x_{0},-x_{1},...,-x_{r}).
\]

\medskip

Suppose now that $G=SO(3)$, then $\Hom(\Z^{n},SO(3))$ is no 
longer connected. In this case $G/T\cong \S^{2}$ with $W$ acting 
as the antipodal map and therefore as in the case of $SU(2)$, 
$G/N(T)\cong \R P^{2}$ and $p_{n}:G/T\times_{W} \t^{n}\to (G/T)/W$ 
is the vector bundle $n\lambda_{2}$. In this case, the map 
$\alpha_{n}$ is a homeomorphism and therefore 
\[
\Hom(\Z^{n},SO(3))_{(1,...,1)}/S^{1}_{n}(SO(3))_{(1,...,1)}
\cong(\R P^{2})^{n\lambda_{2}}.
\]
On the other hand,  by \cite{Torres} 
\begin{align}\label{idenSO(3)}
\Hom(\Z^{n},SO(3))\cong\Hom(\Z^{n},SO(3))_{(1,...,1)}
\sqcup \left(\bigsqcup_{A(n)} \S^{3}/Q_{8}\right),
\end{align}
where
\[
A(n):=\frac{(2^{n}-1)(2^{n-1}-1)}{3}
\] 
and $Q_{8}$ denotes the quaternion group. Moreover, 
if $n\ge 2$ and 
\[
\ux:=(x_{1},...,x_{n})\in \Hom(\Z^{n},SO(3))
\] 
belongs to a path-connected component different from 
$\Hom(\Z^{n},SO(3))_{(1,...,1)}$, with $x_{k}=1$ for some $k$, 
then $y_{k}=1$ for every commuting sequence $\uy:=(y_{1},...,y_{n})$ 
that belongs to the same path-connected component of $\ux$. This 
shows that under the identification (\ref{idenSO(3)}) the subspace 
$S_{n}^{1}(SO(3))$ corresponds to 
\[
S_{n}^{1}(SO(3))\cong S_{n}^{1}(SO(3))_{(1,...,1)}\sqcup 
\left(\bigsqcup_{B(n)} \S^{3}/Q_{8}\right),
\]
where $B(n)$ is the number of connected components homeomorphic 
to $\S^{3}/Q_{8}$ containing an $n$-tuple
$(x_{1},...,x_{n})$ with at least one $x_{i}=1$. It follows then that
\[
\Hom(\Z^{n},SO(3))/S^{1}_{n}(SO(3))\cong 
\Hom(\Z^{n},SO(3))_{(1,...,1)}/S^{1}_{n}(SO(3))_{(1,...,1)}\vee 
\left(\bigvee_{C(n)} (\S^{3}/Q_{8})_{+}\right).
\]
The numbers $C(n)$ are such that $C(1)=0$, $C(2)=1$ 
and satisfy the recurrence relation
\[
\sum_{r=1}^{n}C(r)\binom{{n}}{{r}}=\frac{(2^{n}-1)(2^{n-1}-1)}{3}
\] 
whose solution is 
\[
C(n)=\frac{1}{2}(3^{n-1}-1).
\]
Since $SO(3)\cong \R P^{3}$, it follows that there are homeomorphisms 
\begin{equation*}
\Hom(\Z^{n},SO(3))/S^{1}_{n}(SO(3))\cong \left\{ 
\begin{array}{cl}
\R P^{3} & \text{if }n=1, \\  
(\R P^{2})^{n\lambda_{2}}\vee 
\left(\bigvee_{C(n)} (\S^{3}/Q_{8})_{+}\right)& \text{if }n\ge 2. 
\end{array}%
\right. 
\end{equation*}

Modulo conjugation, the group $SO(3)$ acts transitively on each 
connected component of $\Hom(\Z^{n},SO(3))$ that is homeomorphic to 
$\S^{3}/Q_{8}$. Therefore
\[
\Rep(\Z^{n},SO(3))/\bar{S}^{1}_{n}(SO(3))
\cong(\bigvee_{C(n)} \S^{0})\vee \S^{n}/\Sigma_{2}
\]
with $\Sigma_{2}$ acting on $\S^{n}$ as before. This completely 
describes the factors $\Hom(\Z^{n},G)/S_{n}^{1}(G)$ and 
$\Rep(\Z^{n},G)/\bar{S}_{n}^{1}(G)$ when $G$ is a rank 1 compact, 
connected Lie group. 

\medskip

Next, the case of almost commuting elements for $G$ of rank one 
is discussed. If $G=\S^{1}$ or $SO(3)$ the spaces of almost commuting 
elements and commuting elements agree. Thus the only case left to 
consider is $G=SU(2)$ and $K=Z(SU(2))\cong \Z/2$. The stable factors 
$B_{n}(SU(2),\Z/2)/S^{1}_{n}(SU(2),\Z/2)$ can be described 
in a similar way as in the case of commuting elements. To see this, 
recall that 
\[
B_{n}(SU(2),\Z/2)=\bigsqcup_{C\in T(n,\Z/2)}\AC_{SU(2)}(C).
\]
If $C=C_{1}$ is the trivial matrix whose entries are all equal to $1$, 
then 
\[
\AC_{SU(2)}(C)=\Hom(\Z^{n},SU(2))
\] 
which is a path-connected space. 
On the other hand, suppose that $C=(c_{ij})$ is a non trivial 
matrix with $\AC_{SU(2)}(C)$ not empty. In this case there exist 
$1\le i<j\le n$ such that $c:=c_{ij}\ne 1$. For each $1\le k\le n$ 
write $c_{ik}=c^{b_{k}}$ and $c_{jk}=c^{a_{k}}$ with 
$0\le a_{k},b_{k}\le 1$. In \cite{ACG2} it is 
proved that the map 
\begin{align}\label{map psi}
\psi:\AC_{SU(2)}(C(c))\x (\Z/2)^{n-2}&\to \AC_{SU(2)}(C)\\
((x_{i},x_{j}),(w_{1},...,\hat{w_{i}},...,\hat{w_{j}},...,w_{n}))
&\mapsto(y_{1},...,y_{n}), 
\end{align}
is a homeomorphism, where
\begin{equation*}
y_{k}=\left\{ 
\begin{array}{cl}
x_{i} & \text{if }k=i, \\ 
x_{j} & \text{if }k=j, \\ 
w_{k}x_{i}^{a_{k}}x_{i}^{b_{k}} & \text{if }k\neq i,j,%
\end{array}%
\right. 
\end{equation*}
and $C(c)$ is the $2\times 2$ antisymmetric matrix with entries 
in $\Z/2$ with $c_{12}=c_{21}^{-1}=c$. In addition, if $c\ne 1$, 
then $\AC_{SU(2)}(C(c))\cong PU(2)$. It follows that all 
the path-connected components of $B_{n}(SU(2),\Z/2)$ different 
from $\Hom(\Z^{n},SU(2))$ are homeomorphic to $PU(2)\cong SO(3)$. 
Moreover, from the results in \cite{ACG2} it is easy to see that 
there are 
\[
D(n):=\frac{2^{n-2}(2^{n}-1)(2^{n-1}-1)}{3}
\] 
such connected components. 
Therefore
\begin{equation}\label{almost commuting in SU(2)}
B_{n}(SU(2),\Z/2)\cong \left(\bigsqcup_{D(n)}PU(2)\right)
\sqcup \Hom(\Z^{n},SU(2)).
\end{equation}
On the other hand, if $\ux=(x_{1},...,x_{n})\in \AC_{SU(2)}(C)$ 
is an almost commuting sequence in $SU(2)$ with $x_{k}=1$ for 
some $k$, then either $\ux$ is a commuting sequence in which case 
$\ux\in S_{n}^{1}(SU(2))$ and $C$ is trivial or $C$ is a not trivial. 
If $C$ is not trivial, then there are  $1\le i<j\le n$ such that 
$c:=c_{i,j}\ne 1$. Note that necessarily $i,j\ne k$. Since $x_{k}=1$, 
then $c_{ik}=[x_{i},x_{k}]=1$ and $c_{jk}=[x_{j},x_{k}]=1$, thus 
$a_{k}=b_{k}=0$. It follows that in this case the $k$-th coordinate 
of the map $\psi$ is given by  $w_{k}x_{i}^{a_{k}}x_{i}^{b_{k}}=
w_{k}\in \Z/2$. This shows that if $\uy\in \AC_{SU(2)}(C)$ is in 
the same path-connected component as $\ux$, then $y_{k}=1$. It 
follows that $S^{1}_{n}(SU(2),\Z/2)$ is mapped under the 
homeomorphism in (\ref{almost commuting in SU(2)}) onto the 
subspace 
\[
\left(\bigsqcup_{F(n)}PU(2)\right)\sqcup S_{n}^{1}(SU(2)),
\]
where $F(n)$ is the number of path-connected components of 
$B_{n}(SU(2),\Z/2)$ homeomorphic to $PU(2)$ that contain an 
$n$-tuple $(x_{1},...,x_{n})$ such that $x_{k}=1$ for some $k$. 
Therefore 
\[
B_{n}(SU(2),\Z/2)/S^{1}_{n}(SU(2),\Z/2)\cong 
\Hom(\Z^{n},SU(2))/S_{n}^{1}(SU(2))\vee 
(\bigvee_{K(n)} PU(2)_{+}), 
\]
where the $K(n)$'s 
are integers that satisfy $K(1)=0$, $K(2)=1$ and 
\[
\sum_{r=1}^{n}K(r)\binom{{n}}{{r}}=
\frac{2^{n-2}(2^{n}-1)(2^{n-1}-1)}{3}.
\] 
The solution of this recurrence equation is 
\[
K(n)=\frac{7^{n}}{24}-\frac{3^{n}}{8}+\frac{1}{12}.
\]
Since $PU(2)\cong SO(3)\cong \R P^{3}$ then there are 
homeomorphisms
\begin{equation*}
B_{n}(SU(2),\Z/2)/S^{1}_{n}(SU(2),\Z/2)\cong\left\{
\begin{array}{cl}
\S^{3} & \text{if }n=1, \\  
(\bigvee_{K(n)} \R P^{3}_{+})\vee 
(\R P^{2})^{n\lambda_{2}}/s_{n}(\R P^{2})  & \text{if }n\ge 2. \\  
\end{array}
\right.
\end{equation*}

On the other hand, the situation modulo conjugation simplifies 
substantially.
In this case, the action of $SU(2)$ on the 
connected components of $\bar{B}_{n}(SU(2),\Z/2)$
homeomorphic to $PU(2)$ is transitive. It follows that 
\[
\bar{B}_{n}(SU(2),\Z/2)\cong Z_{n}\sqcup \Rep(\Z^{n},SU(2)),
\]
where $Z_{n}$ is a finite set with 
$D(n)=\frac{2^{n-2}(2^{n}-1)(2^{n-1}-1)}{3}$ points.
This together with Theorem \ref{the case of rep} shows that
\[
\bar{B}_{n}(SU(2),\Z/2)\cong Z_{n}\sqcup (\S^{1})^{n}/\Sigma_{2},
\]
with $\Sigma_{2}$ acting diagonally on 
$(\S^{1})^{n}$ and by complex conjugation on each $\S^{1}$ factor and 
\[
\bar{B}_{n}(SU(2),\Z/2)/\bar{S}^{1}_{n}(SU(2),\Z/2)
\cong(\bigvee_{K(n)} \S^{0})\vee \S^{n}/\Sigma_{2}
\]
where as before,  if $\tau\in\Sigma_{2}-\{1\}$ 
then for a point $(x_{0},...,x_{n})\in \S^{n}$
\[
\tau\cdot (x_{0},...,x_{n})=(x_{0},-x_{1},...,-x_{n}).
\]

\section{Stable factors modulo conjugation}

In this section an explicit description for the stable factors 
$\bar{B}_{r}(G,K)/\bar{S}^{1}_{r}(G,K)$ is provided 
for $r=1,2$ and $3$ and when $G$ is a compact, simple, 
connected and simply connected Lie group. This is achieved using 
the work 
in \cite{BFM}.

\medskip

Take $G$ a compact, connected, simple and simply connected Lie group. 
Fix $T\subset G$ a maximal torus, let $\Phi=\Phi(T,G)$ be the root 
system associated to  $T$, $W$ the Weyl group and fix $\Delta$ a 
set of simple roots for $\Phi$. For each $I\subset \Delta$ there 
is an associated torus $S_{I}$ in $G$ whose lie algebra is 
\[
\t_{I}:=\bigcap_{a\in I}\Ker(a)\subset \t.
\]  
For each such $I$, let $L_{I}=DZ(S_{I})$ be the derived group 
of the centralizer of the torus $S_{I}$ in $G$. Also let 
$F_{I}=S_{I}\cap L_{I}$, this is a finite subgroup of the center 
of $L_{I}$. In \cite{BFM} the spaces $\M_{G}(C)=\AC_{G}(C)/G$ 
are described in the following way. Given $\ux\in\M_{G}(C)$, any 
maximal torus in $Z_{G}(\ux)$ is conjugated to a unique torus of 
the form $S_{I}$. The subset $I=I(\ux)$ only depends on $\ux$ and it 
is a locally constant function of $\ux$. Given $I\subset \Delta$, 
let $\M_{G}^{I}(C)$ be the subspace of conjugacy classes of 
almost commuting $n$-tuples in $G$ of type $C$ whose centralizer 
has a maximal torus conjugate to $S_{I}$. Therefore, each 
$\M_{G}^{I}(C)\subset \M_{G}(C)$ is a union of connected components. 
Moreover, by \cite[Theorem 2.3.1]{BFM} there is a homeomorphism
\begin{equation}\label{BFMgeneral}
(S^{n}_{I}\x_{F_{I}^{n}}\M_{L_{I}}^{0}(C))/W(S_{I},G)
\cong \M_{G}^{I}(C)
\end{equation}
where $\M_{L_{I}}^{0}(C)$ is the moduli space of almost commuting 
$n$-tuples $(x_{1},...,x_{n})$ in $G$ of type $C$ for which 
$Z_{G}(x_{1},...,x_{n})$  has rank $0$ and $W(S_{I},G)$ is the Weyl 
group of the torus $S_{I}$; that is, 
$W(S_{I},G)=N_{G}(S_{I})/Z_{G}(S_{I})$. It follows that for 
$K\subset Z(G)$ a closed subgroup each $\bar{B}_{n}(G,K)$ is a union 
of spaces of the form  (\ref{BFMgeneral}). Moreover, in \cite{BFM} 
explicit description for these spaces that only depend on the geometry 
of $G$ are given for the cases $n=1,2$ and $3$. These descriptions 
can be used to determine the stable factors 
$\bar{B}_{r}(G,K)/\bar{S}^{1}_{r}(G,K)$. To do this note that under 
the assumptions $Z(G)$ is finite and 
thus $\pi_{0}(K)=K$ for all $K\subset Z(G)$. Given $K\subset Z(G)$ a 
subgroup, by (\ref{decomposition for almost}) there is a decomposition 
\[
B_{r}(G,K)=\bigsqcup_{C\in T(r,K)}\AC_{G}(C).
\]
Each $\AC_{G}(C)$ is invariant under the action of $G$, 
thus on the level of orbit spaces
\[
\bar{B}_{r}(G,K)=\bigsqcup_{C\in T(r,K)}\M_{G}(C),
\]
where $\M_{G}(C)=\AC_{G}(C)/G$ is as defined before. Suppose that 
$(x_{1},...,x_{r})\in \AC_{G}(C)$ for some $C\in T(r,K)$. If
$x_{i}=1_{G}$ for some $1\le i\le r$, then $c_{ij}=c_{ji}^{-1}=1_{G}$ 
for all $1\le j\le r$; that is, the elements in the $i$-th column 
and $i$-th row of $C$ are all $1_{G}$. For $0\le i\le r$ define 
$T_{i}(r,K)$  to be the subset of $T(r,K)$ consisting of matrices 
$C$ that have exactly $i$ rows (and hence $i$ columns) whose entries 
are all $1_{G}$. Note that if $C$ is an $r\x r$ antisymmetric matrix 
with $(r-1)$ columns only having $1_{G}$'s, then all the entries of 
$C$ are $1_{G}$ and thus $T_{r-1}(r,K)$ is empty. By definition 
\[
T(r,K)=\bigsqcup_{0\le i\le r}T_{i}(r,K)
\]
and thus 
\[
B_{r}(G,K)=\bigsqcup_{0\le i\le r}
\bigsqcup_{C\in T_{i}(r,K)}\AC_{G}(C).
\]
In general, if $X=\sqcup_{\alpha\in \Lambda}X_{\alpha}$ and 
$A\subset X$ then $X/A=\bigvee_{\alpha\in \Lambda}X_{\alpha}/
(A\cap X_{\alpha})$,
where the usual convention $X/\emptyset=X_{+}$ is used. 
In particular, with the notation 
\[
S_{r}^{1}(G,C):=S^{1}_{r}(G,K)\cap\AC_{G}(C)
\]
there is an identification
\[
B_{r}(G,K)/S^{1}_{r}(G,K)=\bigvee_{0\le i\le r}
\bigvee_{C\in T_{i}(r,K)}\AC_{G}(C)/S_{r}^{1}(G,C).
\]
When $i=r$ there is only one element in $T_{r}(r,\pi_{0}(K))$ 
which is the matrix whose entries are all $1_{G}$ that is 
denoted by $C_{1_{G}}$. In this case $\AC_{G}(C_{1_{G}})=
\Hom(\Z^{n},G)$ and $S_{r}^{1}(G,C_{1_{G}})=S_{r}^{1}(G)$, 
therefore
\[
B_{r}(G,K)/S^{1}_{r}(G,K)=X_{r}(G,K)\bigvee 
\Hom(\Z^{r},G)/S_{r}^{1}(G)
\]
where
\[
X_{r}(G,K)=\bigvee_{0\le i\le r-2}\bigvee_{C\in T_{i}(r,K)}
\AC_{G}(C)/S_{r}^{1}(G,C)
\]
From here it follows that all the stable factors in the 
decomposition of $\Hom(\Z^{n},G)$ appear in the stable 
decomposition $B_{n}(G,K)$ for all closed subgroups 
$K\subset Z(G)$. On the level of orbit spaces
\begin{equation}\label{general case r}
\bar{B}_{r}(G,K)/\bar{S}^{1}_{r}(G,K)=\bar{X}_{r}(G,K)\bigvee
\Rep(\Z^{r},G)/\bar{S}_{r}^{1}(G),
\end{equation}
with 
\[
\bar{X}_{r}(G,K)=\bigvee_{0\le i\le r-2}\bigvee_{C\in T_{i}(r,K)}
\M_{G}(C)/\bar{S}_{r}^{1}(G,C).
\]
The cases $r=1,2$ and $3$ are considered next.

$\bullet$ {\bf{Case r=1}}.
If $r=1$, then $B_{1}(G,K)=G$ and $S^{1}_{1}(G,K)=\{1_{G}\}$, 
thus
\begin{align*}
B_{1}(G,K)/S^{1}_{1}(G,K)&=G,\\
\bar{B}_{1}(G,K)/\bar{S}^{1}_{1}(G,K)&=G/G^{ad}=T/W.
\end{align*}

$\bullet$ {\bf{Case r=2}}.
If $r=2$ then by (\ref{general case r}) 
\[
\bar{B}_{2}(G,K)/\bar{S}^{1}_{2}(G,K)=
\bar{X}_{2}(G,K)\bigvee\Rep(\Z^{2},G)/\bar{S}_{2}^{1}(G),
\]
where
\[
\bar{X}_{2}(G,K)=\bigvee_{C\in T_{0}(K)}
\M_{G}(C)/\bar{S}_{r}^{1}(G,C).
\]
Any matrix $C\in T(K)$ is of the form $C=C(c)$,  where $C(c)$ 
denotes the $2\x 2$ antisymmetric matrix with 
$c_{1,2}=c_{2,1}^{-1}=c\in K$. In particular any $C\in T_{0}(K)$ 
is of the form $C(c)$ for $c\in K-\{1_{G}\}$ and for such a matrix 
$S_{2}^{1}(G,C)=\emptyset$. Therefore 
\[
\bar{X}_{2}(G,K)=\bigvee_{c\in K-\{1_{G}\}}
\M_{G}(C)/\bar{S}_{r}^{1}(G,C)
=\bigvee_{c\in K-\{1_{G}\}}(\M_{G}(C))_{+}.
\]
Each $\M_{G}(C)$ can be described in terms of the geometry 
of $G$ as follows. Consider $\Phi^{\vee}$ the inverse root 
system and $\Delta^{\vee}$ the set of coroots $a^{\vee}$ inverse 
to roots $a\in \Delta$. Each $c\in K\subset Z(G)$ can be written 
in the form $c=\exp(\sum_{a\in \Delta}r_{a}a^{\vee})$,
where each $r_{a}\in \Q$. Given such $c$, following \cite{BFM} 
define $I_{c}=\{a\in \Delta/r_{a}\notin\Z\}\subset \Delta$. 
By \cite[Proposition 4.2.1]{BFM} given $(x,y)\in G^{2}$ with 
$[x,y]=c$, then any maximal torus in $Z_{G}(x,y)$ is conjugate 
to $S_{I_{c}}$. In addition, by \cite[Corollary 4.2.2]{BFM}  
\[
\M_{G}(C(c))\cong (\bar{S}_{I_{c}}\x \bar{S}_{I_{c}})/W(S_{I_{c}},G).
\]
On the other hand, in \cite{Borel} Borel proved that 
for any pair $x$ and $y$ of commuting elements in a simply 
connected compact Lie group, there is a maximal torus that 
contains both $x$ and $y$. Since all maximal tori are conjugate 
it follows that $\Rep(\Z^{2},G)=(T\x T)/W$ and 
\[
\Rep(\Z^{2},G)/\bar{S}^{1}_{2}(G)=(T\wedge T)/W.
\]
Therefore 
\[
\bar{B}_{2}(G,K)/\bar{S}^{1}_{2}(G,K)=
\left(\bigvee_{c\in K-\{1_{G}\}}((\bar{S}_{I_{c}}\x 
\bar{S}_{I_{c}})/W(S_{I_{c}},G))_{+}\right)\vee(T\wedge T)/W
\]

$\bullet$ {\bf{Case r=3}}. The case $r=3$ is similar. 
By (\ref{general case r}) 
\[
\bar{B}_{3}(G,K)/\bar{S}^{1}_{3}(G,K)=
\bar{X}_{3}(G,K)\bigvee\Rep(\Z^{3},G)/\bar{S}_{3}^{1}(G)
\]
with 
\[
\bar{X}_{3}(G,K)=\left(\bigvee_{C\in T_{0}(3,K)}\M_{G}(C)_{+}\right)
\vee \left(\bigvee_{C\in T_{1}(3,K)}\M_{G}(C)/
\bar{S}_{3}^{1}(G,C)\right)
\]
Using the work in \cite{BFM} an explicit description for 
the terms in $\bar{X}_{3}(G,K)$ can be given. It can 
be shown that 
\[
\bar{B}_{3}(G,K)/\bar{S}^{1}_{3}(G,K)=Y_{3,0}\vee 
Y_{3,1}\vee Y_{3,2}\vee (T\wedge T\wedge T)/W
\]
where
\begin{align*}
Y_{3,0}=&\bigvee_{I_{0}}(\bar{S}_{I_{0}}^{3}/W(S_{I_{0}},G))_{+},\\
Y_{3,1}=&\bigvee_{I_{1}}((S^{2}_{I_{1}}\x 
(S_{I_{1}}/K_{I_{1}}))/W(S_{I_{1}},G))_{+},\\
Y_{3,2}=&\bigvee_{I_{2}}(\bar{S}^{2}_{I_{2}}\x 
S_{I_{2}})/(\bar{S}^{2}_{I_{2}}\x\{1_{G}\})/W(S_{I_{2}},G).
\end{align*}
where $I_{0},I_{1}$ and $I_{2}$ run through certain 
subsets of $\Delta$ that are determined uniquely by 
the geometry of $G$ and $K_{I_{1}}$ is a group acting 
on $S_{I}$ of order at most $2$. (In \cite[Theorem 9.5.1]{BFM} 
the different possibilities for $K_{I_{1}}$ are discussed 
depending on the type of $G$).

\medskip

In general, by \cite[Theorem 2.3.1]{BFM} the space 
$\bar{B}_{n}(G,Z(G))$ can be written as a disjoint union of 
spaces of the form 
\begin{equation}\label{BFMgeneral1}
(S^{n}_{I}\x_{F_{I}^{n}}\M_{L_{I}}^{0}(C))/W(S_{I},G).
\end{equation}
It follows that for $n\ge 1$ the stable factor 
$\bar{B}_{n}(G,K)/\bar{S}^{1}_{n}(G,K)$ can be written as the 
wedge of certain quotients of spaces of the form 
(\ref{BFMgeneral1}).

\end{document}